%% file: e6.tex
\documentclass[12pt,a4paper]{article}
\usepackage[american]{babel}
\usepackage{amsmath}
\usepackage{latexsym}
\usepackage{amssymb}
\usepackage{theorem}
\theoremstyle{change}
\newtheorem{Thm}{Theorem}[section]
\newtheorem{Cor}[Thm]{Corollary}
\newtheorem{Prop}[Thm]{Proposition}
\newtheorem{Lem}[Thm]{Lemma}
{\theorembodyfont{\rmfamily}
\newtheorem{Num}[Thm]{}

\newtheorem{Def}[Thm]{Definition}}
\renewcommand{\phi}{\varphi}
\renewcommand{\rho}{\varrho}

\newcommand{\bra}[1]{\langle#1\rangle}
\newcommand{\proof}{\par\medskip\rm\emph{Proof. }}

\newcommand{\qed}{\ \hglue 0pt plus 1filll $\Box$}

\newcommand{\mapstoo}{\longmapsto}

\newcommand{\RR}{\mathbb{R}}

\newcommand{\CC}{\mathbb{C}}
\newcommand{\id}{\mathrm{id}}
\newcommand{\SKIP}[1]{}

\renewcommand{\emptyset}{\varnothing}
\newcommand{\tr}{\mathrm{tr}}

\newcommand{\Aut}{\mathrm{Aut}}

\newcommand{\Spin}{\mathrm{Spin}}

\newcommand{\Inv}{\mathrm{Inv}}

\newcommand{\cP}{\mathcal P}

\newcommand{\cL}{\mathcal L}

\newcommand{\bfi}{\mathbf I}
\newcommand{\calp}{\cP}
\newcommand{\call}{\cL}
\newcommand{\calv}{\mathcal V}
\newcommand{\calk}{\mathcal K}
\newcommand{\calh}{\mathcal H}
\newcommand{\cc}{\mathbb C}
\newcommand{\rr}{\mathbb R}
\newcommand{\hh}{\mathbb H}

\newcommand{\oo}{\mathbb O}
\newcommand{\ff}{\mathbb F}
\newcommand{\re}{\mathrm{Re}}

\newcommand{\conj}[1]{#1^\ast}
\newcommand{\ovl}{\overline}
\newcommand{\ug}{\mathrm U}
\newcommand{\og}{\mathrm O}
\newcommand{\gl}{\mathrm{GL}}

\newcommand{\ext}{\textstyle\bigwedge}
\newcommand{\one}{{(1)}}
\newcommand{\two}{{(2)}}
\newcommand{\herm}{\mathfrak H}
\newcommand{\matr}{\mathfrak M}
\newcommand{\hvec}[2]{(#1_1,#1_2,#1_3;#2_1,#2_2,#2_3)}
\newcommand{\hv}[2]{((#1_i)_i; (#2_i)_i)}
\newcommand{\hermmatr}[2]
  { \left(
    \begin{array}{ccc}
    #1_1 & #2_3 & \conj{#2}_2 \\
    \conj{#2}_3 & #1_2 & #2_1 \\
    #2_2 & \conj{#2}_1 & #1_3
    \end{array}
    \right)
  }
\DeclareMathAlphabet{\varcal}{U}{rsfs}{m}{it}

\begin{document}
\title{{\bf The real quadrangle of type $E_6$}}
\author{Torsten Kurth, Ralf Gramlich, Linus Kramer}
%\date{}
\maketitle

\begin{abstract}
Based on the first author's diploma thesis \cite{Kurth:2000} we use the theories of Lie groups and of Tits buildings in order to describe a Veronese embedding of the real quadrangle of type $E_6$, i.e., the $C_2$~sub-building of the complex $E_6$~building corresponding to the real form $E_{6(-14)}$ of the semisimple complex Lie group of type $E_6$.
\end{abstract}

\section*{Introduction}
The complex Lie group of type $E_6$ has four noncompact real forms,
denoted $E_{6(6)}$, $E_{6(2)}$, $E_{6(-26)}$, and $E_{6(-14)}$
in \cite{Helgason}. The first of these groups is the split
real form of $\RR$-rank $6$, the second has $\RR$-rank $4$ and
belongs to a building of type $F_4$, and the third noncompact real form
is the
automorphism group of the real projective Cayley plane. This group $E_{6(-26)}$
has been studied in detail by Freudenthal \cite{ProjektiveOktavenebene};
a modern account is given in \cite{CompactProjectivePlanes}.

This paper is concerned with the last group, $E_{6(-14)}$.
Veldkamp studied it in his papers on real forms of groups of
type $E_6$ in detail; however, this was in the ``pre-building''
days and was carried out largely in the language of Hjelmslev planes
\cite{HjelmslevI}, \cite{HjelmslevIII}.
The building for the group $E_{6(-14)}$ is a polar space
of type $C_2$, i.e.\ a generalized quadrangle. In an abstract
building-theoretic setting, it is studied in \cite{TitsWeiss}, \cite{WeissE6}, \cite{Weiss:2006}
in detail (over arbitrary fields). It is our aim
to describe a projective embedding (called a Veronese embedding)
of this generalized quadrangle over the reals,
which is similar in spirit to Freudenthal's description
of the Cayley plane. 

Similarly as the Cayley plane, the $E_6$ generalized quadrangle has a
smaller and somewhat less mysterious relative, the dual of a
classical hermitian generalized quadrangle. We first describe
a Veronese embedding for this ``toy model'' before getting
to $E_6$. This dual classical generalized quadrangle exhibits
in fact many of the properties (and difficulties) we encounter
afterwards with $E_6$.

In the first three sections of this article we collect information about incidence geometries of rank two in general and generalized quadrangles in particular (Section  \ref{sec1}), Cayley and Jordan algebras (Section \ref{sec2}), and real and complex Lie groups of types $E_6$ and $F_4$ (Section \ref{sec3}). This information is combined in Section \ref{sec4} in order to describe the Veronese embedding of the real quadrangle of type $E_6$ (Theorem \ref{la.TheFinalTheorem}); in an Appendix \ref{sec5} we provide a concrete list of equations that coordinatize this embedding. The ``toy model'' dual classical hermitian generalized quadrangle mentioned in the introduction is described in \ref{1.15} and \ref{1.16}.   

\medskip \noindent
{\bf Acknowledgements:} This paper is based on the first author's
diploma thesis \cite{Kurth:2000},
written in 2000 in W\"urzburg under
Theo Grundh\"ofer's and the third author's supervision.
In the present form, it was written up by the last two authors.
We are indebted to Arjeh Cohen,
Theo Grundh\"ofer, Bernhard M\"uhlherr and Richard Weiss. We also thank an anonymous referee for helping us to
improve the exposition of the material covered in this article.

\section{Incidence geometries and some linear algebra} \label{sec1}
We collect some basic notions from incidence geometry.
A general reference
is the book \cite{GeneralizedPolygons}.

\begin{Num}\textbf{Incidence geometries}\label{la.IncidenceGeometry}\\
An {\em incidence geometry (of rank $2$)}  is a
triple $\Gamma = (\calp,\call,\bfi)$ consisting of two non-empty disjoint sets
$\calp$ and $\call$ and a reflexive symmetric
relation $$\bfi \subseteq (\calp \cup \call) \times (\calp \cup \call),$$ the {\em incidence relation}, satisfying $\bfi_{| \calp \times \calp} = \mathrm{id}$ and $\bfi_{| \call \times \call} = \mathrm{id}$. 
The elements of $\calp$ and $\call$ are called {\em points}
and {\em lines}, respectively. 
A point $p$ and a line $L$ are {\em incident} if $p \bfi L$. The elements
of $\calv = \calp \cup \call$ are also called {\em vertices}. We say that
$\Gamma$ is {\em thick} if every vertex $x$ of
$\Gamma$ is incident with at least three other vertices.
Two points $p,q\in\calp$ are called {\em collinear} if there exists a line
$L\in\call$ such that $p\bfi L\bfi q$, in which case we write
$p\bot q$.
Similarly we call two lines $L$ and $M$ {\em
confluent}, in symbols $L\bot M$, if there exists a
point $p\in\calp$ satisfying $L\bfi p \bfi M$.
If $\Gamma = (\calp,\call,\bfi)$ is an incidence geometry, then the incidence
geometry $\Gamma^D = (\call,\calp,\bfi)$ 
is called the {\em dual incidence geometry of $\Gamma$}.
\end{Num}
Collineations (isomorphisms) between incidence geometries are
type-preserving and incidence-preserving permutations of the set $\calp \cup \call$ whose inverses are also incidence-preserving. The automorphism group of $\Gamma$
is denoted by $\Aut(\Gamma)$.
\begin{Num}\textbf{Generalized quadrangles}\label{la.GeneralizedQuadrangle}\\
A \emph{projective plane} is a thick geometry where any two distinct points
are joined by a unique line, and any two distinct lines meet in a unique point.
A thick incidence geometry $\Gamma = (\calp,\call,\bfi)$ is called a {\em
generalized quadrangle} if any two distinct lines have at most one
point in common, and if for any non-incident point-line pair
$(p,M)\in\cP\times\cL$ there exists a unique point-line pair
$(q,L)$ such that $p\bfi L\bfi q\bfi M$.

\medskip\noindent
Clearly, the dual incidence geometry of a projective plane
(a generalized quadrangle) is also a
projective plane (a generalized quadrangle).
\end{Num}
For details on generalized quadrangles we refer to \cite{Payne/Thas:1984},
\cite{GeneralizedPolygons}. 

\begin{Num}\textbf{Collinearity relations}
\label{la.CollinearityRelation}\\
Let $S$ be a non-empty set and $R \subseteq S \times S$ a reflexive
and symmetric relation. For $x,y\in S$ let
$L(x,y) = \{ z \in S \mid zRx, zRy \}$. Then
$$\Gamma_R = (S,\{ L(x,y) \mid x,y \in S, x \not= y, xRy \}, \in)$$
is an incidence geometry. If $\Gamma$ is a thick incidence geometry
and if $\Gamma$ contains no digons and triangles, then
obviously $\Gamma_\bot\cong\Gamma$, where $\bot\subseteq\cP\times\cP$
is the collinearity relation of $\Gamma$. This applies in particular
to generalized quadrangles (whereas the collinearity relation of a
projective plane carries no information).
\end{Num}
Assume now that $K$ is a field of characteristic
$\mathrm{char}(K)\neq 2$ and $V$ a vector space over $K$
of finite dimension $n$.
\begin{Num}\textbf{Hermitian forms}\label{la.HermitianForm}\\
Let $\sigma : K \rightarrow K$ be an automorphism of the field $K$ satisfying
$\sigma^2=\id_K$.  A
{\em $\sigma$-hermitian form} is a map
$h:V\times V \rightarrow K$ such that for all $x,y,x_1,x_2,y_1,y_2 \in V$ and
$a,b \in K$ we have
\begin{enumerate}
\item[(i)] $h(ax,by)=ab^\sigma h(x,y)$,
\item[(ii)] $h(x_1+x_2,y_1+y_2) = h(x_1,y_1) + h(x_1,y_2) + h(x_2,y_1) + h(x_2,y_2)$, and
\item[(iii)] $h(x,y) = h(y,x)^\sigma$.
\end{enumerate}

\medskip\noindent
For a subset $M\subseteq V$ we define the linear subspace
$$M^{\bot_h} = M^\bot = \{ y \in V \mid h(x,y) = 0 \text{ for all } x \in M \}.$$
We call $h$ {\em non-degenerate}
if $V^\bot = \{ 0 \}$. A linear subspace $U\subseteq V$
is {\em non-degenerate} if $U^\perp\cap U= \{0 \}$, \emph{degenerate}
if $U^\bot\cap U\neq \{0 \}$, and {\em totally isotropic}
if $U\subseteq U^\perp$.
The maximal dimension of a totally isotropic
subspace of $V$ is called the {\em Witt index}
of $h$.
The \emph{unitary group} of a non-degenerate hermitian form is defined as
$$\ug(V,h) = \{ \varphi\in\gl(V) \mid h(\varphi(x),\varphi(y)) = h(x,y)
\text{ for all } x,y \in V \}$$
\end{Num}
\begin{Num}\textbf{Classical quadrangles}\label{la.ClassicalQuadrangle}\\
Assume that
$n = \mathrm{dim}(V) \geq 5$ and let $h : V \times V \rightarrow K$ be a non-degenerate
$\sigma$-hermitian form of Witt index $2$. Let $\calp$ be the set of all
one-dimensional totally isotropic subspaces of $V$ and $\call$ be the set of
all two-dimensional totally isotropic subspaces of $V$.
We call the incidence geometry
$$Q(V,h) = (\calp, \call, \subseteq)$$
a {\em classical quadrangle (over $K$)}. 
\end{Num}

\begin{Prop} \label{la.ClassicalQuadrangleThm}
The incidence geometry $Q(V,h)$ defined above is a generalized quadrangle.

\proof
See \cite[Chapter~7]{GeometryOfClassicalGroups} or
 \cite[Section 2.3]{GeneralizedPolygons}. \qed
\end{Prop}
Clearly the unitary group $\ug(V,h)$ acts on this generalized quadrangle.
Now we want to describe the dual of this quadrangle, using the Pl\"ucker
embedding.
\begin{Def} \label{la.General}
The {\em Grassmannian of $r$-dimensional subspaces of $V$}
is the set
$$G_r(V) = \{ U \mid U \text{ is a linear subspace of } V, \dim_KU=r \}.$$
The elements of $G_1(V)=P(V)$ are called points and the elements
of $G_2(V)$ are called lines.
\end{Def}
To give an algebraic description of the quadrangle $Q(V,h)^D$ we need some
facts on exterior powers and the exterior algebra. For details
see \cite{Greub} and \cite[Chapter~XIX]{LangAlgebra}.
Recall that the second exterior power $\ext^2V$ of $V$ is a vector space of dimension $n\choose 2$.

\begin{Num}\textbf{The Pl\"ucker embedding}\label{la.PluckerEmbedding}\\
Given linearly independent vectors $x,y \in V$, put
$p(Kx+Ky) = K(x \wedge y)$. This is a well-defined injective map
$p:G_2(V) \rightarrow P(\ext^2V)$, 
see \cite[Chapter~1, Section~5]{PrinciplesOfAlgebraicGeometry}. We set
$\calk = p(G_2(V))\subseteq P(\ext^2V)$. Then one has
$$\textstyle
\calk = \{K(x \wedge y) \mid x,y \in V \text{ are linearly independent} \} =
\{ Ku\in P(\ext^2V) \mid u \wedge u = 0 \}$$
(cf. \cite[Chapter~1, Section~5]{PrinciplesOfAlgebraicGeometry} or
\cite[Chapter~1, Section~4.1]{BasicAlgebraicGeometryI}). Note that the second
equality depends on our assumption that $\mathrm{char}(K) \not= 2$.
\end{Num}

\begin{Lem} \label{la.ConfluenceInProjectiveSpace}
Suppose $L_1,L_2 \in G_2(V)$ and $p(L_i) = Ku_i$ where $i=1,2$.
Then $L_1 \cap L_2 \not= \{0\}$ holds if and only if $u_1 \wedge u_2 = 0$.

\proof
For $i=1,2$, let $u_i = x_i \wedge y_i$ where $x_i$ and $y_i$ is a basis
for the
two-dimensional subspace $L_i$. Then we have
\begin{eqnarray*}
  L_1 \cap L_2 \not= \{0\} & \Longleftrightarrow & x_1, y_1, x_2, y_2
    \text{ are linearly dependent} \\
  & \Longleftrightarrow & x_1 \wedge y_1 \wedge x_2 \wedge y_2 = 0 \\
  & \Longleftrightarrow & u_1 \wedge u_2 = 0
\end{eqnarray*}
which proves the lemma.
\qed
\end{Lem}
Suppose that $h:V\times V \rightarrow K$ is a $\sigma$-hermitian
form. Then there exists a unique $\sigma$-hermitian
form $h_2: \ext^2V \times \ext^2V \rightarrow K$ satisfying
$$h_2(x_1 \wedge x_2, y_1 \wedge y_2) = \det((h(x_i,y_j)_{i,j=1,2})) =
h(x_1,y_1)h(x_2,y_2) - h(x_1,y_2)h(x_2,y_1)$$
for all $x_1,x_2,y_1,y_2 \in V$ (see \cite[5.9]{Greub}).
We remark that if $h$ is non-degenerate on
$V$, then $h_2$ is also non-degenerate.
%The next lemma shows that in general the Pl\"ucker embedding does not map
%totally isotropic subspaces in $G_2(V)$ onto totally
%isotropic subspaces contained in $P(\ext^2V)$.
\begin{Lem}
Assume $Ku \in \calk$ and let $h:V \times V \rightarrow K$ be a
$\sigma$-hermitian form. Then one has $h_2(u,u) = 0$ if and only if 
the subspace $p^{-1}(Ku) \in G_2(V)$ is degenerate.

\proof
Let $u = x \wedge y$ such that $x$ and $y$ span $L = p^{-1}(Ku)$. Then
$$
  h_2(u,u) = \det \left(
  \begin{array}{cc}
    h(x,x) & h(x,y) \\ h(y,x) & h(y,y)
  \end{array}
  \right) = h(x,x)h(y,y) - h(x,y)h(y,x)
$$
is a scalar multiple of the discriminant of $L$ (with respect to $h$)
(compare \cite[p.~299]{ClassicalGroups}). Now the lemma follows
from \cite[6.1.9~(i)]{ClassicalGroups}.
\qed
\end{Lem}
Therefore we introduce the following terminology.
\begin{Def} \label{la.WeakStrongIsotropyI}
Let $Ku \in \calk$ and let $h:V\times V\rightarrow K$ be an arbitrary non-degenerate
$\sigma$-hermitian form. Then we call $Ku$
\begin{enumerate}
\item[(i)] {\em weakly isotropic} if $h_2(u,u)=0$ and
\item[(ii)] {\em strongly isotropic} 
if $p^{-1}(Ku)$ is totally isotropic.
\end{enumerate}
\end{Def}
Note that a strongly isotropic subspace is also weakly isotropic.
%Since we assume that $h$ is hermitian and $\mathrm{char}(K)\neq2$,
%a subspace $L\in G_2(V)$ is totally isotropic if and only if
%it is isotropic, i.e. if and only if $h(x,x)=0$ for all $x\in L$.
%This allows us to prove the following lemma.
\begin{Lem}
A space $Ku\in\calk$ is strongly isotropic if and only if for each $Kv\in\calk$ satisfying $u \wedge v = 0$ we have $h_2(u,v) = 0$.
\proof
Assume $p^{-1}(Ku)$ is totally isotropic, let $x_1$, $x_2$ be a basis of $p^{-1}(Ku)$, and let $y_1$, $y_2$ be a basis of $p^{-1}(Kv)$. Since $u \wedge v = 0$, we can choose $x_2 = y_1$ by Lemma \ref{la.ConfluenceInProjectiveSpace}. Then we have $h(x_1,y_1) = h(x_1,x_2) = 0$ and $h(x_2,y_1) = h(x_2,x_2) = 0$. Hence $h_2(x_1 \wedge x_2,y_1 \wedge y_2) = h(x_1,y_1)h(x_2,y_2) - h(x_1,y_2)h(x_2,y_1) = 0$. 

If $p^{-1}(Ku)$ is not totally isotropic, there exists 
$Kv \in \calk \backslash \{ Ku \}$ with $u \wedge v = 0$ such that 
$p^{-1}(Kv) \cap p^{-1}(Ku)$ is not contained in the radical of 
$p^{-1}(Ku)$ and such that $p^{-1}(Kv)^{\perp_h} \cap p^{-1}(Ku) = \{ 0 \}$. 
Choose $0 \neq x_1 \in p^{-1}(Ku)$ with 
$x_1 \not\perp_h p^{-1}(Kv) \cap p^{-1}(Ku)$, choose 
$0 \neq x_2 = y_1 \in p^{-1}(Kv) \cap p^{-1}(Ku)$ and choose 
$0 \neq y_2 \in x_1^{\perp_h} \cap p^{-1}(Kv)$. Then $h(x_1,y_2) = 0$ 
and $h(x_1,y_1) \neq 0 \neq h(x_2,y_2)$. Therefore 
$h_2(x_1 \wedge x_2,y_1 \wedge y_2) = h(x_1,y_1)h(x_2,y_2) - 
h(x_1,y_2)h(x_2,y_1) \neq 0$. 
\qed
\end{Lem}
%
%
%***************************************************
%
%
%Therefore a subspace $L\in G_2(V)$ is totally isotropic if and
%only if every subspace $L'\in G_2(V)$ which intersects $L$
%non-trivially is degenerate *** that's a wrong statement; $L'$ has to lie in $L^{\bot_h}$ *** consequently, \ref{stronglyisotropic} and \ref{1.16} are wrong ***.
%\begin{Lem} \label{stronglyisotropic}
%A weakly isotropic space $Ku\in\calk$ is strongly
%isotropic if and only if every
%subspace $Kv\in\calk$ with $u\wedge v=0$ is weakly isotropic.
%\qed
%\end{Lem}
\begin{Num} \label{la.DescriptionOfDualClassicalQuadrangle}
We are now in the position to give a description of dual classical quadrangles
using suitable exterior algebras. So let $Q(V,h)^D$ be such a quadrangle where
$h:V\times V\rightarrow K$ is a non-degenerate $\sigma$-hermitian form of Witt
index $2$. Then put
$$\calh = \{ Ku\in\calk \mid Ku \text{ is strongly isotropic} \}.$$
The elements of $\calh$ obviously correspond to the points of $Q^D(V,h)$.
Define on $\calh$ a collinearity relation $\bot$ by
$$Ku \bot Kv \Longleftrightarrow u \wedge v = 0.$$
By \ref{la.ConfluenceInProjectiveSpace} the relation $\bot$ is isomorphic to
the collinearity relation in $Q^D(V,h)$, so
$\Gamma_\bot$ is isomorphic to $Q(V,h)^D$, that is,
we have derived the desired description.
\end{Num}
\begin{Num} \label{la.SpecialCase}
Now we specialize to the case where $K=\cc$ and
$V=\cc^6$. Let $\sigma$ denote complex conjugation on $\cc$ and define a non-degenerate $\sigma$-hermitian form $h$ of
Witt index $2$ on $V$ and a non-degenerate symmetric bilinear form $(\cdot,\cdot)$ by setting
\begin{eqnarray*}
h(x,y) & = & -x_1\ovl{y_1} - x_2\ovl{y_2} + x_3\ovl{y_3} + x_4\ovl{y_4} +
  x_5\ovl{y_4} + x_5\ovl{y_5} + x_6\ovl{y_6} \\
(x,y) & = & x_1y_1 + x_2y_2 + x_3y_3 + x_4y_4 + x_5y_5 + x_6y_6
\end{eqnarray*}
for all $x=(x_1,\ldots,x_6), y=(y_1,\ldots,y_6) \in \cc^6$. Clearly, $h$ and
$(\cdot, \cdot)$ are non-degenerate on $\cc^6$. Hence their induced forms
on  $\ext^2\cc^6$ are also non-degenerate.
Note that $h_2$ (cf.\ the paragraph after Lemma \ref{la.ConfluenceInProjectiveSpace}) has Witt index $7$ on the $15$-dimensional space
$\ext^2\cc^6$.
\end{Num}
\begin{Num} \label{1.15}
We now give a description of the quadrangle $Q(\cc^6,h)^D$ in the vector space
$\ext^2\cc^6$ using a cross-product. Fix a vector space isomorphism
$\varphi :\ext^6\cc^6 \rightarrow \cc$ and define a trilinear
form $(\cdot,\cdot,\cdot)$ on $\ext^2\cc^6$ by
$$(u,v,w) = \varphi(u \wedge v \wedge w)$$
for all $u,v,w\in\ext^2\cc^6$. Note that this form is symmetric
(i.e. invariant under all permutations of the variables) and
non-degenerate.
Now we can define a bilinear and symmetric
cross-product $\times$ on $\ext^2\cc^6$ by
$$(u \times v,w) = (u,v,w) \quad \text{for all }w \in \ext^2\cc^6$$
for $u,v \in \ext^2\cc^6$.
Note also that
$$u \times v = 0 \Longleftrightarrow u \wedge v = 0$$
for all $u,v \in \ext^2\cc^6$.
\end{Num}
\begin{Num} \label{1.16}
To summarize, we have
\begin{eqnarray*}
\calk&=&\{Ku\in P(\ext^2V)\mid u\times u=0\}\\
\calh&=&\{Ku\in\calk\mid h_2(u,v)=0 \text{ for all }
Kv\in\calk\text{ with }u\times v=0\}\\
%\calh&=&\{Ku\in\calk\mid h_2(u,u)=0=h_2(v,v)\text{ for all }
%Kv\in\calk\text{ with }v\times u=0\}\\
\bot&=&\{(Ku,Kv)\in\calh\times\calh\mid u\times v=0\}.
\end{eqnarray*}
So we have derived a description of the dual classical quadrangle
$Q(\cc^6,h)^D$ in the $15$-dimensional vector space
$\ext^2\cc^6$ in terms of a symmetric trilinear form,
a symmetric bilinear form, a cross-product (defined by these forms)
and a $\sigma$-hermitian form $h_2$. Note also that
$\ext^2V$ is a $\ug(V,h)$-module, so $\ug(V,h)$ acts on
$Q(\cc^6,h)^D$.
\end{Num}

\section{Cayley and Jordan algebras} \label{sec2}

One of the purposes of this section is to remind the reader of the constructions of the algebras of real and complex
quaternions and octonions by using the so-called Cayley-Dickson process. 
A general reference for the material covered below are
e.g.\ \cite{BasicAlgebraI}, \cite{LangAlgebra}, \cite{CompactProjectivePlanes}. Later in this section we introduce Jordan algebras which we will then use to give a description of the $E_6$ building and the Cayley plane. 

We continue to assume that
$K$ is a field of characteristic $\mathrm{char}(K)\neq 2$.

\begin{Num}\textbf{Algebras}\label{la.Algebra}\\
An {\em algebra over $K$} or {\em $K$-algebra} \index{algebra} is a
$K$-vector space $A$ equipped with a $K$-bilinear map
$(x,y) \mapsto  x \cdot y = xy$.
An element $1_A \in A$ is called
an {\em identity element of $A$} \index{identity element} if
$1_Ax = x1_A = x$ holds for all $x \in A$. If
$xy=1_A=yx$, then $x,y$ are called invertible. If the maps 
$x\mapsto ax$ and $x\mapsto xa$ are bijective for every
$a\in A\setminus\{0\}$, then $A$ is
called a \emph{division algebra}. 
An algebra is called \emph{alternative} if it satisfies
the following weak form of associativity for all $x,y\in A$
\[
x^2y=x(xy)\text{ and }yx^2=(yx)x.
\]
A map $\phi\in\gl(A)$ is called an \emph{automorphism} if
$\phi(1_A)=1_A$ and if 
$\phi(xy)=\phi(x)\phi(y)$ holds for all $x,y\in A$, and an
\emph{anti-automorphism}
if instead
$\phi(xy)=\phi(y)\phi(x)$ holds for all $x,y\in A$.
An anti-automorphism whose square is the identity is also
called an \emph{involution}.
\end{Num}
%\begin{Lem} \label{la.AlternativityFormula}
%Let $A$ be an alternative $K$-algebra. Then the
%following equations hold for all $x,y,z\in A$
%
%(1) $x(yz) + y(xz) = (xy)z + (yx)z$
%
%(2) $x(yz) + x(zy) = (xy)z + (xz)y$.
%
%\proof
%Expand both sides of
%$(x+y) ((x+y)z) = (x+y)^2z$ and of
%$x(y+z)^2 = (x(y+z))(y+z)$.
%\qed
%\end{Lem}
Recall that the real Cayley algebra $\oo$ is an $8$-dimensional
alternative division algebra over the reals (in fact, it is the
unique non-associative alternative division algebra over $\rr$).
We briefly recall the construction of $\oo$.
\begin{Num} \label{la.CayleyDickson}
The Cayley-Dickson process yields a family of $K$-algebras
$\ff_m^K$ with unit and with a canonical involutions
$x \mapsto \conj{x}$. We begin with
$\ff_0^K = K$, $x \mapsto \conj{x} := x$.
Assuming that the $K$-algebra $\ff_m^K$ is defined for $0 \leq m$,
we put
$\ff_{m+1}^K = \ff_m^K \oplus \ff_m^K$ with the product
$$(x_1,x_2) \cdot (y_1,y_2) = (x_1y_1 - \conj{y_2}x_2, x_1y_2 +
x_2\conj{y_1}).$$ The involution $x \mapsto \conj{x}$ extends 
to $\ff_{m+1}^K$ via
$$\conj{{(x_1,x_2)}} = (\conj{x_1}, -x_2).$$
Clearly, $\dim_K\ff_m^K=2^m$. Via the embedding
$x\mapsto(x,0)$ we can view $\ff_m^K$ as a subalgebra of
$\ff_{m+1}^K$, the involutions commute with this embedding.
In this way, we have embeddings
$K=\ff_0^K\subseteq\ff_1^K\subseteq\cdots\subseteq\ff_m^K$ and
$1\in K$ is a unit element. Moreover, $K=\ff_0^K$ is a central
subalgebra in $\ff_m^K$ and consists precisely of the fixed elements
of the involution $x\mapsto\conj{x}$.
\end{Num}
\begin{Num}
As usual, we put
$$\re_K(x) = \frac{1}{2}(x+\conj{x})$$
for $x \in \ff_m^K$. Note that $\re_K(x)=\re_K(\conj x) \in K$ 
and that $\re_\rr$ commutes with the inclusions $\ff_m\subseteq\ff_{m+1}$.
We also define the \emph{Norm form} by
$$N_K(x) = x\conj{x}$$ and its polarization
\[
\bra{x|y}=N_K(x+y)-N_K(x)-N_K(y)=2\re_K(x\conj y).
\]
Note that every element $x$ with $N_K(x)\neq 0$ is invertible,
since $x\conj x=\conj xx=N_K(x)$. For
$z=(x,y)\in\ff_{m+1}^K$ we have $N_K(z)=N_K(x)+N_K(y)$.
\end{Num}
Obviously the Cayley-Dickson process is functorial, an inclusion of
fields $K\subseteq L$ yields algebra inclusions
$\ff_m^K\subseteq\ff_m^L$, and there is a natural
isomorphism $\ff_m^K\otimes_KL=\ff_m^L$.
In particular, the complex conjugation $z\mapstoo\ovl z$ on $\cc$
extends to an automorphism of $\ff_m^\cc$ (which we denote by the same
symbol).
\begin{Num} \label{la.Complexification}
We continue with some more observations. Since
$x^2=2\re_K(x)x-N_K(x)$, we see that every element is
contained in an associative and commutative subalgebra.
Also, $\ff_1^K$ is obviously commutative. Direct inspection
shows: if $\ff_m^K$ is commutative, then $\ff_{m+1}^K$ is
associative, and if $\ff_m^K$ is associative, then
$\ff_{m+1}^K$ is alternative. So we see that
$\ff_2^K$ is associative and $\ff_3^K$ is alternative.
The algebras $\ff_m^K$ for $m\geq 4$ are of no further interest to us.
Note that the algebra $\ff_2^K$ is in fact the quaternion
algebra associated to the quaternion symbol $\left(\frac{-1,-1}K\right)$.
The non-associative alternative algebra $\ff_3^K$ is
a Cayley algebra over $K$. Any two elements $x,y\in\ff_3^K$
are contained in an associative subalgebra; in particular,
$\ff_3^K$ is a division algebra if and only if the quadratic
form $N_K$ is anisotropic.
This is the case for $K=\rr$; the real Cayley division algebra is
\[
\oo=\ff_3^\rr.
\]
Note that $\ff_1^\rr=\cc$, and $\ff_2^\rr$ is the (unique) real
quaternion division algebra.
For $K=\cc$, already $\ff_1^\cc\cong\cc[\varepsilon]/(\varepsilon^2+1)$
has zero divisors. We put
$\hh^\cc=\ff_2^\cc=\hh\otimes_\rr\cc$ and
$\oo^\cc=\ff_3^\cc=\oo\otimes_\rr\cc$. This is the complex
Cayley algebra; as we noted, it has zero divisors.

The following identities hold for all $x,y,z \in \ff_m^K$ for $m\leq 3$.
\begin{itemize} 
\item $N_K(xy) = N_K(x)N_K(y)$,
\item $\bra{xy | z} = \bra{ x | \conj{y}z}$,
\item $\bra{xy | z} = \bra{ yx | z }$, 
$\bra{x | yz} = \bra{ x | zy }$,
\item $\re_K(xy) = \re_K(yx)$,
\item $\re_K(x(yz)) = \re_K((xy)z)$, so we may write $\re_K(xyz)$.
\end{itemize}
\end{Num}

\begin{Num}
Next we introduce the Jordan algebra $\herm_3^K$
of hermitian $(3\times 3)$-matrices with entries in the Cayley algebra.
Then on this  Jordan algebra we introduce a
symmetric bilinear form and a symmetric trilinear form which we use to define
a cross-product on $\herm_3^K$ and the notion of a Veronese vector (see
\ref{la.ABilinearForm}, \ref{la.ATrilinearForm}, \ref{la.TheCrossProduct} and
\ref{la.VeroneseVector}). This cross-product and the Veronese vectors turn out
to be important ingredients in the description of the so-called $E_6$ building
over $\cc$ (compare \ref{la.RemarkOnE6Building}). The section closes with a
short description of the real projective octonion plane.
\end{Num}
From now on we assume that $\mathrm{char}(K)\neq2,3$.

\begin{Num} \label{la.Matrix}
Consider the matrix algebra $\matr_n(K)=K^{n\times n}$ of
$(n\times n)$-matrices over $K$. The matrix transposition 
$X\mapsto X^T$ is an involution
on this algebra. Given any $K$-algebra $A$, let
$\matr_n(A)=A\otimes_K\matr_n(K)$. The elements of this algebra
are $(n\times n)$-matrices with entries in $A$, with the usual
matrix multiplication. If $A$ admits an involution $x\mapsto\conj x$,
then $a\otimes X\mapsto\conj a \otimes X^T$ is an involution on
$\matr_n(A)$ which we denote by the same symbol.

We apply this remark with $A = \ff_m^K$. Let
\[
\herm_n^K=\{X\in\matr_n(\ff_3^K)\mid X=\conj X\}.
\]
As a $K$-vector space, $\dim\herm_n^K=(4n-3)n$.
We define a product $\circ$ on this vector space
by
\[
X \circ Y = \frac{1}{2}(XY + YX).
\]
Obviously, this makes $\herm_n^K$ into a commutative algebra
with unit element $I$ (the $n\times n$ identity matrix).
\end{Num}
\begin{Thm} \label{la.JordanAlgebra}
The $K$-vector space $\herm_3^K$ equipped with the product $\circ$
is a commutative 27-dimensional Jordan algebra over $K$
(and hence power-associative).

\proof
See \cite[Chapter 6, \S 4 and Chapter 7, \S 6]{JordanAlgebren}.
\qed
\end{Thm}
We write a typical element $X \in \herm_3^K$ as
$$X = \hermmatr{\xi}{x}$$
or, shorter, as
$$X = \hvec{\xi}{x} = \hv{\xi}{x}$$
where $\xi_i \in K$ and $x_i \in \ff_3^K$ for $i=1,2,3$.
\begin{Prop} \label{la.JordanProduct}
For $X=\hv{\xi}{x}$ and $Y=\hv{\eta}{y}$,
the product $X \circ Y=Z=\hv{\zeta}{z}$ is given by the formula
\begin{eqnarray*}
  \zeta_i & = & \xi_i\eta_i + \frac{1}{2} \langle x_j|y_j \rangle + \frac{1}{2}
    \langle x_k|y_k \rangle \\
  z_i & = & \frac{1}{2} ( (\xi_j+\xi_k)y_i + (\eta_j+\eta_k)x_i +
    \conj{{(y_jx_k)}} + \conj{{(x_jy_k)}} )
\end{eqnarray*}
for $(i,j,k) \in \{(1,2,3),(2,3,1),(3,1,2)\}$.

\proof
Write
$$X = \hermmatr{\xi}{x} \quad \text{and} \quad Y=\hermmatr{\eta}{y}$$
and compute
$$XY = \left(
  \begin{array}{ccc}
  \xi_1\eta_1 + x_3\conj{y}_3 + \conj{x}_2y_2 &
    \conj{x}_1y_3 + x_3\eta_2 + \conj{x}_2\conj{y}_1 &
    \xi_1\conj{y}_2 + x_3y_1 + \conj{x}_2\eta_3 \\
  \conj{x}_3\eta_1 + \xi_2\conj{y}_3+x_1y_2 &
    \conj{x}_3y_3 + \xi_2\eta_2 + x_1\conj{y}_1 &
    \conj{x}_3\conj{y}_2 + \xi_2y_1 + x_1\eta_3 \\
  x_2\eta_1 + \conj{x}_1\conj{y}_3 + \xi_3y_2 &
    x_2y_3 + \conj{x}_1\eta_2 + \xi_3\conj{y}_1 &
    x_2\conj{y}_2+\conj{x}_1y_1+\xi_3\eta_3
  \end{array}
  \right) .
$$
The results follows, since $Z=\frac{1}{2}(XY+YX)$.
\qed
\end{Prop}
\begin{Cor} \label{la.JordanProductSquare}
For $X=\hv{\xi}{x}$, we have
$X \circ X=XX=Z=\hv{\zeta}{z}$, with
\begin{eqnarray*}
  \zeta_i & = & \xi_i^2 + N_K(x_j) + N_K(x_k) \\
  z_i & = & (\xi_j+\xi_k)x_i + \conj{{(x_jx_k)}}
\end{eqnarray*}
for $(i,j,k) \in \{ (1,2,3), (2,3,1), (3,1,2) \}$. \qed
\end{Cor}
\begin{Num}\textbf{The bilinear form}\label{la.ABilinearForm}\\
For $X,Y,Z\in \herm_3^K$ put
$\tr(Z)=\zeta_1+\zeta_2+\zeta_3$ and
$$(X,Y) = \tr(X \circ Y).$$ Note that $(X,I)=\tr(X)$.
This is obviously a symmetric
bilinear form on the vector space $\herm_3^K$. 
From \ref{la.JordanProduct} we get
$$(X,Y) = \xi_1\eta_1 + \xi_2\eta_2 + \xi_3\eta_3 + \langle x_1|y_1 \rangle +
\langle x_2|y_2 \rangle + \langle x_3|y_3 \rangle,$$
hence $(\cdot,\cdot)$ is non-degenerate. It can be shown
that $(X,Y \circ Z) = (X \circ Y, Z)$ (see \cite[Chapter 7, \S 5]{JordanAlgebren}), so we may write
$$\tr(X \circ (Y \circ Z)) = \tr((X \circ Y) \circ Z)
= \tr(X \circ Y \circ Z).$$
\end{Num}
\begin{Num}\textbf{The trilinear form}\label{la.ATrilinearForm}\\
For $X,Y,Z \in \herm_3^K$ we define a symmetric trilinear form $(\cdot,\cdot,\cdot) : \herm_3^K \times \herm_3^K \times \herm_3^K \to K$ by
\begin{eqnarray*}
  3(X,Y,Z) & = & \tr(X \circ Y \circ Z) - \frac{1}{2} \tr(X) \tr(Y \circ Z) - \frac{1}{2}
    \tr(Y)\tr(X \circ Z) \\
  & & - \frac{1}{2} \tr(Z)\tr(X \circ Y) + \frac{1}{2} \tr(X)\tr(Y)\tr(Z).
\end{eqnarray*}
Furthermore, the {\em determinant of $X$} is defined as
$$\det X = (X,X,X) = \frac{1}{3} \tr(X \circ X \circ X) - \frac{1}{2}
\tr(X)\tr(X^2) + \frac{1}{6} {(\tr(X))}^3.$$
One can check that
$$6(X,Y,Z) = \det(X+Y+Z) - \det(X+Y) - \det(X+Z) - \det(Y+Z) + \det X +
\det Y + \det Z.$$
\end{Num}
\begin{Num}\textbf{The cross product}\label{la.TheCrossProduct}\\
We define a symmetric bilinear map
$\times:\herm_3^K\times\herm_3^K\to\herm_3^K$
through
$$(X \times Y,Z) = 3 (X,Y,Z).$$
Then
\begin{eqnarray*}
  (X \times Y,Z) & = & \tr(X \circ Y \circ Z) - \frac{1}{2} \tr(X) \tr(Y \circ Z) - \frac{1}{2} 
    \tr(Y)\tr(X \circ Z) \\
  & & - \frac{1}{2} \tr(Z)\tr(X \circ Y) + \frac{1}{2} \tr(X)\tr(Y)\tr(Z) \\
  & = & \tr\left( \left( X \circ Y - \frac{1}{2} \tr(X)Y - \frac{1}{2}\tr(Y)X - \frac{1}{2} \tr(X
    \circ Y)I \right. \right. \\
  & & \qquad \left. \left. + \frac{1}{2} \tr(X)\tr(Y)I \right) \circ Z \right) \\
  & = & \left( X \circ Y - \frac{1}{2} \tr(X)Y - \frac{1}{2}\tr(Y)X - \frac{1}{2}
    \tr(X \circ Y)I \right. \\
  & & \qquad \left. + \frac{1}{2} \tr(X)\tr(Y)I, Z \right)
\end{eqnarray*}
whence
$$X \times Y = X \circ Y - \frac{1}{2} (Y,I)X - \frac{1}{2} (X,I)Y -
\frac{1}{2} (X,Y)I + \frac{1}{2} (X,I)(Y,I)I.$$
\end{Num}
\begin{Prop} \label{la.CrossProduct}
For $Z = X \times Y$ we have
\begin{eqnarray*}
  2\zeta_i & = & \xi_j\eta_k + \xi_k\eta_j - \langle x_i|y_i \rangle \\ 
  2z_i & = & \conj{{(y_jx_k)}} + \conj{{(x_jy_k)}} - \xi_iy_i - \eta_ix_i
\end{eqnarray*}
where $(i,j,k) \in \{ (1,2,3),(2,3,1),(3,1,2) \}$.

\proof
A direct computation using the last formula.
\qed
\end{Prop}
\begin{Cor} \label{la.CrossProductSquare}
For $Z = X \times X$ we have
\begin{eqnarray*}
  \zeta_i & = & \xi_j\xi_k - N_K(x_i) \\ 
  z_i & = & \conj{{(x_jx_k)}} - \xi_ix_i
\end{eqnarray*}
where $(i,j,k) \in \{ (1,2,3),(2,3,1),(3,1,2) \}$.
\end{Cor}
\begin{Cor} \label{la.Determinant}
We have
$$\det X = \xi_1\xi_2\xi_3 - \xi_1N_K(x_1) - \xi_2N_K(x_2) - \xi_3N_K(x_3) +
2 \re(x_1x_2x_3).$$

\proof
Expand
$3 \det X =  (X \times X, X)$.
\qed

\end{Cor}
\begin{Cor} 
We have
$$X\det X =(X\times X)\times(X\times X).$$

\proof
Expand both sides.
\qed
\end{Cor}
\begin{Def} \label{la.VeroneseVector}
A non-zero element $X=\hv{\xi}{x} \in \herm_3^K$ is
called a {\em Veronese vector} if $X \times X = 0$.
Note that this implies $\det(X)=0$ (see \ref{la.ATrilinearForm} and \ref{la.TheCrossProduct}).
By \ref{la.CrossProductSquare} an element $X\neq0$
is a Veronese vector if and only
if the {\em Veronese conditions}
\begin{eqnarray*}
  N_K(x_i) & = & \xi_j\xi_k, \\
  \xi_i\conj{x_i} & = & x_jx_k
\end{eqnarray*}
are satisfied for all $(i,j,k) \in \{(1,2,3), (2,3,1), (3,1,2)\}$. Finally,
we define the set
$$\calv = \{ K X \mid X \text{ is a Veronese vector} \}$$
of all one-dimensional subspaces of $\herm_3^K$ which are generated by
Veronese vectors.
\end{Def}

\begin{Num}\textbf{The building of type $E_6$}\label{la.RemarkOnE6Building}\\
Suppose now that the $8$-dimensional
norm form $N_K$ is isotropic over $K$
(for example that $K$ is quadratically or algebraically closed, e.g.~$K=\cc$).
Then $\ff_3^K$ is the split Cayley algebra, and
it follows that the product $\times$ given in \ref{la.TheCrossProduct}
is identical with the cross-product defined in
\cite[p.~691pp]{PointLineSpaces}; this follows from the results proved in
\cite{ClassOfJordanAlgebras}, \cite{ClassOfCubicForms}.
Then the set $\calv$ given in
\ref{la.VeroneseVector} is the point set of an incidence geometry, the
{\em $E_6$ building over $\cc$} which is
defined in \cite[p.~691pp]{PointLineSpaces}; see also \ref{la.TheE6Building} of the present paper. Moreover, two points
$KX, K Y \in \calv$ of the $E_6$ building are {\em collinear}
\index{collinear} if and only if $X \times Y = 0$ (compare also
\cite[p.~692]{PointLineSpaces}).
\end{Num}

\begin{Num}\textbf{The Cayley plane}\label{la.ProjectiveOctonionPlane}\\
On the other hand, assume that the norm form on
$\ff_3^K$ is anisotropic; this holds for example if
$K$ is formally real, eg.~$K=\rr$. Then $\ff_3^K$ is
a Cayley division algebra. We define
\begin{eqnarray*}
  \calp & = & \calv,  \\
  \call & = & \{ X^{\bot_{(\cdot,\cdot)}} \mid KX \in \calp \}.
\end{eqnarray*}
For $K=\rr$, the incidence geometry
$$\Gamma = (\calp, \call, \subseteq)$$
is the so-called (and well-known) projective Cayley plane which is
studied in \cite[Sections 16, 17, 18]{CompactProjectivePlanes} and
Freudenthal's classical paper \cite{ProjektiveOktavenebene}, and
Adams' book \cite{Adams}.
\end{Num}

\section{Some exceptional Lie groups} \label{sec3}

In this section we define several noteworthy groups of automorphisms of the real and complex Jordan algebras defined in the preceding section.

\begin{Num} {\bf The invariance group and automorphisms}\label{la.InvarianceGroup}\\ \label{la.AutomorphismGroup}
Let $K \in \{\rr,\cc\}$. We consider the Jordan algebra $\herm_3^K$ with its bilinear form $(\cdot,\cdot)$ (see \ref{la.ABilinearForm}) and determinant $\mathrm{det}$ (see \ref{la.ATrilinearForm}). The
set
\begin{eqnarray*}
\Inv_K(\det) & = & \{ \varphi\in\gl(\herm_3^K) \mid \det \varphi(X) = \det X
\text{ for all } X \in \herm_3^K \} \\
& = & \{ \varphi\in\gl(\herm_3^K) \mid (\varphi(X),\varphi(Y),\varphi(Z)) = 
      (X,Y,Z) \text{ for all } X,Y,Z \in \herm_3^K \}
\end{eqnarray*}
forms a subgroup of $\gl(\herm_3^K)$. (Compare \ref{la.ATrilinearForm} for the latter equality.) Clearly, $\Inv_\rr(\det) \leq \Inv_\cc(\det)$. 

By \cite{ProjektiveOktavenebene},
\cite[Section~17]{CompactProjectivePlanes} the group $\Inv_\rr(\det)$ is a simple real Lie group of type
$E_{6(-26)}$. 
The group $E_{6(-26)}$ is the collineation group of the real projective
Cayley plane described in \ref{la.ProjectiveOctonionPlane} (see
\cite{CompactProjectivePlanes}, Section~17). Moreover, by \cite[Section~7]{GeometricAlgebraProjectiveOctavePlanes} the group $\Inv_\cc(\det)$ is an almost simple complex Lie group of type $E_6$.

It can be shown (cf.\ \cite{ProjektiveOktavenebene}) that
$\Aut(\herm_3^K) = \Inv_K(\det) \cap \og(\herm_3^K, (\cdot,\cdot))$. Since $\Inv_\rr(\det) \leq \Inv_\cc(\det)$ we have $\Aut(\herm_3^\rr) \leq \Aut(\herm_3^\cc)$. By \cite{ProjektiveOktavenebene} and
\cite[Section~18]{CompactProjectivePlanes} the group $\Aut(\herm_3^\rr)$ is a simple real Lie group of
type $F_{4(-52)}$. By \cite[Section~7]{GeometricAlgebraProjectiveOctavePlanes} the group $\Aut(\herm_3^\cc)$ is a simple complex Lie group of type $F_4$.
\end{Num}

\begin{Prop}
Let $\varphi \in \Inv_\cc(\det)$ and $X,Y \in \herm_3^\cc$. Then $X \times Y = 0$ if and only if $\varphi(X) \times \varphi(Y) = 0$.
\proof
Since $\varphi$ lets the trilinear form
$(\cdot,\cdot,\cdot)$ invariant (see \ref{la.InvarianceGroup}), we have
\begin{eqnarray*}
  (\varphi(X) \times \varphi(Y), T) & = & 3 (\varphi(X),\varphi(Y),T) = 
    3 (\varphi(X),\varphi(Y),\varphi(\varphi^{-1}(T))) \\
  & = & 3 (X,Y,\varphi^{-1}(T)) = (X \times Y, \varphi^{-1}(T))
\end{eqnarray*}
for all $T \in \herm_3^\cc$. Because $T$ is arbitrary and $\phi$ is a bijection, the non-degeneracy of
$(\cdot,\cdot)$ implies the desired equivalence.
\qed
\end{Prop}

\begin{Num} {\bf Real forms of $E_6$ and $F_4$}\label{la.ARealFormOfE6}\\
Recall that $z\mapstoo\ovl z$ denotes the involutive automorphism of the complex
Cayley algebra $\oo_\CC=\oo\otimes\cc$ induced by complex conjugation on the scalars.
The map
$$H : \herm_3^\cc \rightarrow \herm_3^\cc : \hvec{\xi}{x} \mapsto (\ovl{\xi}_1, \ovl{\xi}_2, \ovl{\xi}_3; - \ovl{x}_1,
- \ovl{x}_2, \ovl{x}_3)$$ is $\rr$-linear,
$\cc$-semilinear, and bijective. Furthermore, $H$
preserves the Veronese conditions, so that it maps Veronese vectors onto
Veronese vectors. Moreover, given $X,Y \in \herm_3^\cc$, one has
$X \times Y = 0$ if and only if $H(X) \times H(Y) = 0$. 

The form $h : \herm_3^\cc \times \herm_3^\cc \to \cc$ defined by
$$h(X,Y) = (X, H(Y)) = \xi_1\ovl{\eta}_1 + \xi_2\ovl{\eta}_2 +
\xi_3\ovl{\eta}_3 - \langle x_1|\ovl{y}_1 \rangle - \langle x_2|\ovl{y}_2
\rangle + \langle x_3|\ovl{y}_3 \rangle$$
for all $X=\hv{\xi}{x},Y=\hv{\eta}{y} \in \herm_3^\cc$ is $\bar\ $-hermitian. By \cite[Sections 5, 6, 7]{HjelmslevIII} the group
$\Inv_\cc(\det) \cap \ug(\herm_3^\cc,h)$ is an almost simple real Lie group of type $E_6$ (in fact $E_{6(-14)}$ by Theorem \ref{e614} below). Moreover, by \cite[Section~18]{CompactProjectivePlanes} the group $\Inv_\rr(\det) \cap \og(\herm_3^\rr,h_{|\herm_3^\rr \times \herm_3^\rr})$ is a simple real Lie group of type $F_{4(-20)}$.
\end{Num}

\begin{Lem} \label{la.IdentificationOfRealForm}
An almost simple real Lie group of type $E_6$ which
contains an almost simple real Lie group of type $F_{4(-20)}$ and the group
$\Spin(10)$ as subgroups, is of type $E_{6(-14)}$.
\proof
This follows by inspection of the ranks and dimensions of the real forms of the complex Lie group of type $E_6$.
\qed
\end{Lem}

\begin{Num} \label{la.ARealFormOfF4}
We turn our attention to $\herm_3^\rr$. We define an automorphism
$B \in \gl(\herm_3^\rr)$ by restricting $H$ to this real subspace,
$$B(\hvec{\xi}{x}) = (\xi_1,\xi_2,\xi_3; -x_1,-x_2,x_3)$$
for all $\hv{\xi}{x} \in \herm_3^\rr$ and a symmetric bilinear form
$\beta$ on $\herm_3^\rr$ by putting
$$\beta(X,Y) = (X,B(Y)) = \xi_1\eta_1 + \xi_2\eta_2 + \xi_3\eta_3 - \langle
x_1|y_1 \rangle - \langle x_2|y_2 \rangle + \langle x_3|y_3 \rangle$$
for all $X=\hv{\xi}{x},Y=\hv{\eta}{y} \in \herm_3^\rr$.
Then we have
$$H(X \otimes \xi) = H((\xi^\one X, \xi^\two X)) =  B(X) \otimes \ovl{\xi}$$
and therefore
$$h(X \otimes \xi, Y \otimes \eta)= \beta(X,Y) \xi\ovl{\eta},$$ where $h$ is the form defined in \ref{la.ARealFormOfE6}.
From this we obtain $\og(\herm_3^\rr,\beta) \leq \ug(\herm_3^\cc,h)$.
Hence we can consider the group $\Inv_\rr(\det) \cap \og(\herm_3^\rr,\beta)$
as a subgroup of $\Inv_\cc(\det)\cap\ug(\herm_3^\cc,h)$.
\end{Num}

\begin{Num} \label{la.SpinGroup}
{\bf The groups $\Spin(8)$, $\Spin(9)$, $\Spin(10)$} Let $a \in \oo$ with $N_\RR(a)=1$. We define a map $T_a \in \gl(\herm_3^\cc)$
by
$$T_a(X) =
  \left(
    \begin{array}{ccc}
      a & 0 & 0 \\ 0 & \conj{a} & 0 \\ 0 & 0 & 1
    \end{array}
    \right)
    \hermmatr{\xi}{x}
    \left(
    \begin{array}{ccc}
      \conj{a} & 0 & 0 \\ 0 & a & 0 \\ 0 & 0 & 1
    \end{array}
    \right)
  =
    \left(
    \begin{array}{ccc}
      \xi_1 & ax_3a & \conj{{(x_2\conj{a})}} \\
      \conj{{(ax_3a)}} & \xi_2 & \conj{a}x_1 \\
      x_2\conj{a} & \conj{{(\conj{a}x_1)}} & \xi_3
    \end{array}
    \right)
  $$
for all $X=\hv{\xi}{x} \in \herm_3^\cc$. The group generated by the maps
$T_a$ is the group $\Spin(8)$ (see \cite[p.~267]{CliffordAlgebras}).
Furthermore, for an ordered pair $(c,s) \in \rr^2$ where $c^2+s^2=1$ one
defines a map $R_{(c,s)} \in \gl(\herm_3^\cc)$ by setting
\begin{eqnarray*}
R_{(c,s)}(X) & = &
    \left(
    \begin{array}{ccc}
      c & s & 0 \\ -s & c & 0 \\ 0 & 0 & 1
    \end{array}
    \right)
    \hermmatr{\xi}{x}
    \left(
    \begin{array}{ccc}
      c & -s & 0 \\ s & c & 0 \\ 0 & 0 & 1
    \end{array}
    \right) \\
 & = &
    \left(
    \begin{array}{ccc}
      \begin{array}{c}
        c^2\xi_1 + s^2\xi_2 \\
        + 2cs\re(x_3) 
      \end{array}
      &
      \begin{array}{c}
        c^2x_2 - s^2\conj{x_3} \\
        + cs(\xi_2-\xi_1)
      \end{array}
      & \conj{{(cx_2 + s\conj{x_1})}} \\
      \hfill \\
      \begin{array}{c}
        \conj{{(c^2x_2 - s^2\conj{x_3})}} \\
        + cs(\xi_2-\xi_1)
      \end{array}
      &
      \begin{array}{c}
        s^2\xi_1 + c^2\xi_2 \\
        - 2cs\re(x_3)
      \end{array}
      & -s\conj{x_2} + cx_1 \\
      \hfill \\
      cx_2 + s\conj{x_1} & \conj{{(-s\conj{x_2} + cx_1)}} & \xi_3
    \end{array}
    \right)
\end{eqnarray*}
for every $X=\hv{\xi}{x} \in \herm_3^\cc$. The group generated by the maps
$T_a$ and $R_{(c,s)}$ is the group $\Spin(9)$ (compare \cite[p.~275]{CliffordAlgebras}). It is well known that $\Spin(9)$ is a subgroup of $F_{4(-52)}$ and
hence of the complex Lie group $\mathrm{Aut}(\herm_3^\cc) = \Inv_\cc(\det) \cap \og(\herm_3^\cc,(\cdot,\cdot))$
(see \ref{la.InvarianceGroup}). The maps $T_a$ and
$R_{(c,s)}$ obviously commute with the map $H$ defined in
\ref{la.ARealFormOfE6}, that is, we have $H \circ T_a = T_a \circ H$ and
$H \circ R_{(c,s)} = R_{[c,s)} \circ H$. Because of $\Spin(9) \leq
\og(\herm_3^\cc, (\cdot,\cdot))$ this implies $\Spin(9) \leq
\ug(\herm_3^\cc,h)$. Thus $\Spin(9)$ is a subgroup of $\Inv_\cc(\det) \cap \ug(\herm_3^\cc,h)$. A similar
argument shows by restricting every $\varphi \in \Spin(9)$ to $\herm_3^\rr$
and seeing that $B$ (cf.\ \ref{la.ARealFormOfF4}) commutes with each generator of $\Spin(9)$ that
$\Spin(9)$ is a subgroup of $F_{4(-20)}$. \par
We now consider maps $S_\omega \in \gl(\herm_3^\cc)$ where
$\omega \in \cc$ satisfies $\omega \ovl{\omega} = 1$ which are defined by
$$S_\omega(X) =
    \left(
    \begin{array}{ccc}
      \omega & 0 & 0 \\ 0 & \ovl{\omega} & 0 \\ 0 & 0 & 1
    \end{array}
    \right)
    \hermmatr{\xi}{x}
    \left(
    \begin{array}{ccc}
      \omega & 0 & 0 \\ 0 & \ovl{\omega} & 0 \\ 0 & 0 & 1
    \end{array}
    \right)
  =
    \left(
    \begin{array}{ccc}
      \omega^2\xi_1 & x_3 & \omega\conj{x_2} \\
      \conj{x_3} & \ovl{\omega}^2\xi_2 & \ovl{\omega}x_1 \\
      \omega x_2 & \ovl{\omega}\conj{x_1} & \xi_3
    \end{array}
    \right)
$$
for all $X=\hv{\xi}{x} \in \herm_3^\cc$. The group generated by
$\Spin(9)$ and the maps $S_\omega$ is the group $\Spin(10)$ (see
\cite[p.~282f]{CliffordAlgebras}). From \ref{la.Determinant} and
\ref{la.ARealFormOfE6} one directly sees that each map $S_\omega$ is
an element of $\Inv_\cc(\det) \cap \ug(\herm_3^\cc,h)$. Hence
$\Spin(10)$ (and also $F_{4(-20)}$) is a subgroup of $\Inv_\cc(\det) \cap \ug(\herm_3^\cc,h)$. Thus \ref{la.IdentificationOfRealForm} allows us to determine the type of the group $\Inv_\cc(\det) \cap \ug(\herm_3^\cc,h)$.
\end{Num}

\begin{Thm} \label{e614}
The group $\Inv_\cc(\det) \cap \ug(\herm_3^\cc,h)$ is an almost
simple Lie group over the reals of type $E_{6(-14)}$. \qed
\end{Thm}

\section{The real quadrangle of type $E_6$} \label{sec4}

Recall from \ref{la.VeroneseVector} that $\calv = \{ \cc X \mid X \text{ is a Veronese vector} \}$.
As remarked in \ref{la.RemarkOnE6Building}, the set $\calv$ is the point set of the
complex building of type $E_6$. In \ref{la.TheE6Building} we will describe an incidence geometry related to this building in detail.

\begin{Prop}
The group $\Inv_\cc(\det)$ acts transitively on the set
$\calv$, as does any co-central quotient.
\proof
See \cite[Section 3]{HjelmslevI}.
\qed
\end{Prop}

\begin{Num}{\bf Isotropic points}\label{la.IsotropicPoint}\\
Let $p = \cc X \in \calv$. Recall the definition of $H$ and $h$ from  \ref{la.ARealFormOfE6}. We call $p$ {\em weakly isotropic} if $h(X,X) = 0$ and {\em strongly isotropic} if additionally $4 H(X) \times (X \times T) = h(T,X)X$
for all $T \in \herm_3^\cc$ (compare \cite[Definition~(1.3)]{HjelmslevIII}, \cite[Section~1]{HjelmslevI}). Denote the set of strongly isotropic points by $\calh$.
\end{Num}

\begin{Thm}
The real Lie group $\Inv_\cc(\det) \cap \ug(\herm_3^\cc,h)$ of type $E_{6(-14)}$ acts transitively on the set $\calh$
of all  strongly isotropic points.
\proof
See \cite[Section~5]{HjelmslevIII}.
\qed
\end{Thm}

\begin{Num} \label{la.E6QuadrangleI}
Let $p,q \in \calh$ be strongly isotropic points. Write $p = \cc X$ and
$q = \cc Y$ where $X$ and $Y$ are Veronese vectors. Put
$p \bot q :\Longleftrightarrow X \times Y = 0$ (see \ref{la.TheCrossProduct}).
Clearly, this defines a collinearity relation $\bot$ on $\calh$. The aim of
this section is to prove that the incidence geometry $\Gamma_\bot$
associated to $\bot \subset \calh \times \calh$ is a generalized quadrangle, namely the real
$E_6$ quadrangle $Q(E_6,\rr)$. To this end we prove that the geometry $\Gamma_\bot$ is isomorphic to the sub-building of type $C_2$ fixed by the involution $\iota$ defined in \ref{inv} of the complex $E_6$ building; see \ref{la.CollinearityE6Quadrangle} and \ref{la.TheFinalTheorem}.
\end{Num}

\begin{Num}{\bf Singular subspaces of $P(\herm_3^\cc)$}\label{la.SingularSubspace}\\
Let $S = P(U)$ be a subspace of $P(\herm_3^\cc)$ where $U$ is a linear
subspace of $\herm_3^\cc$. Call $S$ {\em singular} provided that $X \times Y = 0$ holds for all
$X,Y \in U$. Note that the vectors in $U$ are Veronese vectors. A singular
subspace $S$ is called {\em maximal} \index{subspace!maximal singular} if
for every singular subspace $S'$ satisfying $S \subseteq S'$ we have $S=S'$.
It can be shown that each maximal singular subspace of $P(\herm_3^\cc)$
has rank $4$ or $5$, and both cases occur; see \cite[p.~259p]{HjelmslevI}
and \cite[p.~693]{PointLineSpaces}.
\end{Num}

\begin{Num}{\bf Symplecta}\label{la.Symplecton}\\
Suppose that $X$ is a Veronese vector. Then
$$X \times \herm_3^\cc = \{ X \times Y \mid Y \in \herm_3^\cc \}$$
is a linear subspace of $\herm_3^\cc$. Call a subset $Q$ of $\calv$ a
{\em symplecton} \index{symplecton} provided that $Q$ is of the form
$P(X \times \herm_3^\cc) \cap \calv$ where $X$ is a Veronese vector.
Note that the symplecton $Q$ generates the projective space $P(X \times
\herm_3^\cc)$, and, vice versa, every projective space of the form
$P(X \times \herm_3^\cc)$, where $X$ is a Veronese vector, determines a
unique symplecton $Q$ which is given by $Q = P(X \times \herm_3^\cc)
\cap \calv$. So in the sequel we shall not distinguish between symplecta and
such projective spaces. If we put
$$\call_1 = \calv = \{ \cc X \mid X \text{ is a Veronese vector} \} \quad
\text{and} \quad \call_6 = \{ Q \mid Q \text{ is a symplecton} \},$$
then
$$\cc X \mapsto P(X \times \herm_3^\cc) \cap \calv$$
defines a bijective map from $\call_1$ onto $\call_6$ such that the following
holds: If one has $\cc X, \cc Y \in \call_1$ where $X$ and $Y$ are Veronese
vectors, then $X \times Y = 0$ holds (that is, $\cc X$ and $\cc Y$ are
collinear) if and only if $P(X \times \herm_3^\cc) \cap \calv$ and
$P(Y \times \herm_3^\cc) \cap \calv$ meet in a maximal singular subspace
of rank $4$ (cf. \cite[p.~694]{PointLineSpaces}).
The intersection of two distinct symplecta is either empty or a single point
or a maximal singular subspace (in each symplecton); see
\cite[p.~694]{PointLineSpaces}.
\end{Num}

\begin{Num}{\bf The complex $E_6$ building}\label{la.TheE6Building}\\
We define an incidence geometry $\Gamma$ of rank 6. For this, we put
\begin{eqnarray*}
  \call_3 & = & \{ l \mid l \text{ is a singular subspace of rank } 1 \}, \\
  \call_4 & = & \{ E \mid  E \text{ is a singular subspace of rank } 2 \}, \\
  \call_5 & = & \{ y \mid y \text{ is a maximal singular subspace of rank } 4 \}, \\
  \call_2 & = & \{ v \mid v \text{ is a maximal singular subspace of rank } 5 \} \\
\end{eqnarray*}
and define incidence relations $\bfi_{ij} \subseteq \call_i \times \call_j$
where $1 \leq i < j \leq 6$ in the following way. Let $x \in \call_i$ and
$y \in \call_j$. If $(i,j) \not= (2,5), (2,6)$, define
$x \bfi_{ij} y$ if and only if $x \in y$.
Moreover,
$x \bfi_{25} y$ if and only if $x$ and $y$ meet in a
singular subspace of rank $3$,
and
$x \bfi_{26} y$ if and only if $x$ and $y$ meet in a
non-maximal singular subspace of rank $4$.
We usually suppress the indices and write $x \bfi y$ instead of $x \bfi_{ij} y$.
Then
$$\Gamma = (\call_1, \call_2, \call_3, \call_4, \call_5, \call_6,
( \bfi_{ij} )_{1 \leq i < j \leq 6} )$$
is an incidence geometry of rank six, whose chamber system is a complex spherical
building of type $E_6$
(see \cite[p.~696]{PointLineSpaces}), i.e., a building related to the diagram
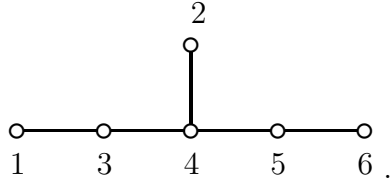
\begin{center} \input{e6i}. \end{center}
The complex Lie group of type $E_6$ acts as a chamber-transitive group of automorphisms on this building.
As usual we
call the elements of $\call_3$ {\em lines}. \index{line} We recall that two
points $\cc X$ and $\cc Y$ are
collinear if and only if $X \times Y = 0$ is
satisfied (compare \ref{la.RemarkOnE6Building}). Moreover, we define for a
point $p = \cc X \in \call_1$ the subspace
$p^\bot = P( \{ Y \in \herm_3^\cc \mid X \times Y = 0 \} )$.
Finally, we remark that for $2 \leq i \leq 6$ each element $x \in \call_i$
is uniquely determined by the set of all points incident with $x$.
\end{Num}

\begin{Num} {\bf The involution $\iota$}\label{inv}\\
We consider the map $\cc X \mapsto P(X \times \herm_3^\cc) \cap
\calv$ from $\call_1$ onto $\call_6$ as described in \ref{la.Symplecton}. This
map can be extended to an involution $\iota_0$ of the building $\Gamma$ such that
$\iota_0$ satisfies the following conditions (see \cite{DiagramGeometries}, 
\cite[5.3]{PointLineSpaces}):
\begin{enumerate}
\item $\iota_0(\call_2) = \call_2$, $\iota_0(\call_4) = \call_4$.
\item $\iota_0(\call_1) = \call_6$, $\iota_0(\call_6) = \call_1$.
\item $\iota_0(\call_3) = \call_5$, $\iota_0(\call_5) = \call_3$.
\end{enumerate}
The function $H$ from \ref{la.ARealFormOfE6} maps $\call_1$ bijectively onto itself and
preserves collinearity. Since $\Gamma$ is
determined by its point set and the collinearity of points (compare \ref{la.TheE6Building}), this map $H$
induces a collineation of $\Gamma$, denoted by $\varphi_H$. Then
$\iota = \varphi_H \circ \iota_0$ is a permutation of $\Gamma$ satisfying the
above properties; since $\varphi_H$ is an involution and since $\varphi_H$ and $\iota_0$ commute, the permutation $\iota$ is in fact an involution.
By \cite[Sections 5, 6]{HjelmslevIII} the centralizer of the restriction $\iota_{16}$ of $\iota$ to $\call_1 \cup \call_6$ in the complex Lie group $\Inv_\cc(\det)$ of type $E_6$ is a real form of type $E_{6(-14)}$. Since every element
$x \in \call_i$, $2 \leq i \leq 5$, is uniquely determined by the set of
points incident with $x$, this real form also centralizes $\iota$.
By \cite[p.~534, p.~518]{Helgason} this means that $\iota$ is an involution which is related to the diagram
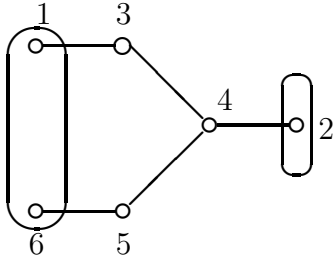
\begin{center} \input{e6ii} .\end{center}  
(We prefer to use the M\"uhlherr diagram instead of the Satake diagram because of our application of results from \cite{MHerr} in \ref{la.E6QuadrangleII} below.)
Consequently, the involution $\iota$ additionally satisfies the following conditions:
\begin{enumerate}
\setcounter{enumi}{3}
\item There exist elements $x \in \call_2$ such that $\iota(x) = x$ but no
elements of $\call_4$ having this property.
\item There exist $\iota$-invariant flags in $\bfi_{16}$ but no flags in
$\bfi_{35}$ having this property.
\end{enumerate}
Here a flag
$(x,y) \in \bfi_{ij}$ is called {\em $\iota$-invariant}
if $\iota(x,y) = (x,y)$ holds.
\end{Num}

\begin{Num}{\bf The real $E_6$ quadrangle}\label{la.E6QuadrangleII}\\
Let $\iota$ be the involution discussed in \ref{inv}. Define
\begin{eqnarray*}
  \calp & = & \{ p \in \call_1 \mid (p, \iota(p)) \text{ is a } \iota
    \text{-invariant flag} \}, \\
  \call & = & \{ L \in \call_2 \mid \iota(L) = L \}.
\end{eqnarray*}
By \cite[Theorems 1.7.27 and 1.8.22]{MHerr} the incidence geometry $Q(E_6,\rr) = (\calp,\call,\bfi_{12})$ is a spherical building of type $C_2$. (We refer the reader to \cite[Chapter 2]{MHerr} for details on how to determine the type of a fixed building using M\"uhlherr diagrams.) In other words, the geometry $Q(E_6,\rr) = (\calp,\call,\bfi_{12})$ is the
{\em real quadrangle of type $E_6$}. Unfortunately, \cite{MHerr} is not easily accessible; an alternative reference dealing with fixed buildings is \cite[Proposition 14.6.1]{DiagramGeometries}, which at the time of writing of this article at least can be accessed via the internet.
\end{Num}
We will now prepare the proof that the geometry $\Gamma_\bot$ defined in \ref{la.E6QuadrangleI} is in fact isomorphic to the real $E_6$ quadrangle $Q(E_6,\rr) = (\calp,\call,\bfi_{12})$.

\begin{Prop} \label{la.PropertyI}
Let $p \in \call_1$ and $P \in \call_6$. Assume that $p$ is not incident with
$P$. Then either $p^\bot \cap P = \emptyset$ or there exists a
unique $v \in \call_2$ such that $p \bfi v \bfi P$. This unique $v$ is
generated by $p$ and the subspace $p^\bot \cap P$.
\proof
Assume $p^\bot \cap P \not= \emptyset$. Then the subspace
$p^\bot \cap P$ is a singular subspace of rank $4$ (see
\cite[p.~694]{PointLineSpaces}). Hence the subspace $v$ generated by $p$ and
$p^\bot \cap P$ is a maximal singular subspace of rank $5$. Clearly, $p$ is
contained in $v$ and $v \cap P = p^\bot \cap P$ holds. Hence we obtain
$p \bfi v \bfi P$. The uniqueness of $v$ is shown in \cite[p.~581]{PolarSpacesOfRank3}.
\qed
\end{Prop}

\begin{Prop} \label{la.PropertyII}
Let $\iota$ be the involution defined in \ref{inv}
and suppose that $(p,P), (q,Q) \in \bfi_{16}$ are $\iota$-invariant flags where
$p \not= q$. Then $p$ is not incident with $Q$, and $q$ is not incident with
$P$.
\proof
We assume that $p \bfi Q$ holds. Applying $\iota$ yields $\iota(Q) = q \bfi P =
\iota(p)$. Hence $p,q \in P \cap Q$. Since $P \cap Q$ is a (maximal) singular
subspace (compare \ref{la.Symplecton}), $p$ and $q$ must be collinear. Let
$l$ be the line joining $p$ and $q$ and put $y = \iota(l)$. Because of
$p,q \bfi l$ we have $y \bfi P,Q$, and therefore we obtain $y = P \cap Q$
which yields $l \bfi y = \iota(l)$. Hence $(l,y) \in \bfi_{35}$ is a
$\iota$-invariant flag which contradicts the properties of $\iota$ established in \ref{inv}. This
completes the proof.
\qed
\end{Prop}

\begin{Thm} \label{la.CollinearityE6Quadrangle}
Let $\iota$ be the involution defined in \ref{inv} and assume
that $p$ and $q$ are points of $Q(E_6,\rr)$. Then $p$ and $q$ are collinear
in $Q(E_6,\rr)$ if and only if $p$ and $q$ are collinear in $\Gamma$.
\proof
"$\Longrightarrow$": Assume that $p$ and $q$ are collinear in $Q(E_6,\rr)$.
Then there exists a $v \in \call_2$ such that $p,q \bfi v$. Hence $p,q \in v$.
Since $v$ is a maximal singular subspace, it follows that $p$ and $q$ are
collinear in $\Gamma$. \par
"$\Longleftarrow$": Let $p$ and $q$ be collinear in $\Gamma$. Put
$P = \iota(p)$ and $Q = \iota(q)$. Then $(p,P)$ and $(q,Q)$ are
$\iota$-invariant flags. Suppose that $p \not= q$ (in the case $p=q$ the
assertion is clear). From \ref{la.PropertyII} we get that $p$ is not incident
with $Q$ and $q$ is not incident with $P$. Since $q \in p^\bot \cap Q$ holds,
there exists a unique $v \in \call_2$ such that $p \bfi v \bfi Q$ is satisfied
(see \ref{la.PropertyI}). Similarly, one obtains a unique $v' \in \call_2$
where $q \bfi v' \bfi P$. Applying of $\iota$ yields $\iota(Q) = q \bfi
\iota(v) \bfi P = \iota(p)$ and hence $\iota(v) = v'$. The theorem is proved
if $v=v'$ is showed. For this, we remark that $P$ can be considered as a polar
space if one takes $\call_1 \cap P$ as point set, $\call_3 \cap P$ as line set
and $\bfi_{13}$ as incidence relation (because the residue of a symplecton is
a diagram geometry of type $D_5$ and hence a polar space). Since $p$ and $q$
are collinear, $P \cap Q$ must be a maximal singular subspace of rank $4$ in
$P$. Thus $p^\bot \cap P \cap Q$ is a singular subspace of rank $3$. Because
of $p \not\in Q$ the singular subspace generated by $p$ and
$p^\bot \cap P \cap Q$ has rank $4$. Since $v$ is generated by $p$ and
$p^\bot \cap Q$, this subspace is the intersection of $v$ and $P$ which implies
$v \bfi P$. Moreover, from $q \in p^\bot \cap Q$ it follows $q \bfi v$. Hence
we have $q \bfi v \bfi P$ and $q \bfi v' \bfi P$, and the uniqueness of
$v'$ yields $v=v'$.
\qed
\end{Thm}

\begin{Num} \label{la.Veldkamp}
Let $X$ and $Y$ be Veronese vectors. We call the point $\cc X$ and the
symplecton $P(Y \times \herm_3^\cc) \cap \calv$ {\em V-incident},
\index{V-incident} in symbols $\cc X \bfi_V P(Y \times \herm_3^\cc)$,
provided that $(X,Y) = 0$ and
$$4 Y \times (X \times T) = (T, Y) X$$
holds for all $T \in \herm_3^\cc$ (see \cite[Sections
1 and 3]{HjelmslevI}). Note that a point $\cc X$ is strongly isotropic if and only if
$\cc X$ is V-incident with $P(H(X) \times \herm_3^\cc) \cap \calv$.
The complex Lie group of type $E_6$ preserves the relation $\bfi_V$ (cf. \cite[Sections 2, 3]{HjelmslevI}).

\medskip \noindent
{\bf Claim.}
$\cc X \bfi_V P(Y \times \herm_3^\cc) \cap \calv \Longleftrightarrow
\cc X \bfi_{16} P(Y \times \herm_3^\cc) \cap \calv$

\begin{proof}
Consider the Veronese vector $U = (1,0,0; 0,0,0)$ and let $X=\hv{\xi}{x},
Z=\hv{\zeta}{z}$. Assume $Z = X \times U$. Then we compute using
\ref{la.CrossProduct}
\begin{eqnarray*}
  \zeta_1 & = & \xi_2 \cdot 0 + \xi_3 \cdot 0 - \langle 0 | 0 \rangle = 0, \\
  \zeta_2 & = & \xi_3 \cdot 1 + \xi_1 \cdot 0 - \langle 0 | 0 \rangle = \xi_3, \\
  \zeta_3 & = & \xi_1 \cdot 0 + \xi_2 \cdot 1 - \langle 0 | 0 \rangle = \xi_2, \\
  z_1 & = & \conj{{(x_2 \cdot 0)}} + \conj{{(0 \cdot x_3)}} - \xi_1 \cdot 0 - 1 \cdot x_1
    = -x_1, \\
  z_2 & = & \conj{{(x_3 \cdot 0)}} + \conj{{(0 \cdot x_1)}} - \xi_2 \cdot 0 - 0 \cdot x_2
    = 0, \\
  z_3 & = & \conj{{(x_1 \cdot 0)}} + \conj{{(0 \cdot x_2)}} - \xi_3 \cdot 0 - 0 \cdot x_3
    = 0.
\end{eqnarray*}
Hence we obtain
$$U \times \herm_3^\cc = \{ (0, \xi_2, \xi_3; x_1, 0, 0) | \xi_2,\xi_3
\in \cc, x_1 \in \oo^\cc \}.$$
Moreover, we have
$$
  X \circ U = \frac{1}{2} (XU + UX) = \frac{1}{2} \left( \left(
  \begin{array}{ccc}
    \xi_1 & 0 & 0 \\ \conj{x_3} & 0 & 0 \\ x_2 & 0 & 0
  \end{array}
  \right) + \left(
  \begin{array}{ccc}
    \xi_1 & x_3 & \conj{x_2} \\ 0 & 0 & 0 \\ 0 & 0 & 0
  \end{array}
  \right) \right) = \frac{1}{2} \left(
  \begin{array}{ccc}
    2 \xi_1 & x_3 & \conj{x_2} \\ \conj{x_3} & 0 & 0 \\ x_2 & 0 & 0
  \end{array}
  \right) .
$$
Hence the set of all $X \in \herm_3^\cc$ satisfying $X \circ U = 0$
equals $U \times \herm_3^\cc$. From \cite[Proposition (1.3) (ii)]{HjelmslevI} it follows that a point $\cc X$ is V-incident with
$P(U \times \herm_3^\cc) \cap \calv$ if and only if $X \circ U = 0$
holds. Hence $\cc X \bfi_V P(U \times \herm_3^\cc) \cap \calv$ is
equivalent to $\cc X \bfi_{16} P(U \times \herm_3^\cc) \cap \calv$. \par
Now assume that $X$ and $Y$ are arbitrary Veronese vectors. Put $p = \cc X$,
$P = P(Y \times \herm_3^\cc) \cap \calv$ and $Q = P(U \times
\herm_3^\cc) \cap \calv$. Assume that $p \bfi_V P$ holds. Since
the complex Lie group of type $E_6$ acts transitively on $\call_6$, it contains an element
$\varphi$ such that $P^\varphi = Q$. Thus $p^\varphi \bfi_V Q$,
because the complex Lie group of type $E_6$ preserves the relation $\bfi_V$ (see \ref{la.Veldkamp}). Hence we have
$p^\varphi \bfi_{16} Q$ and therefore $p \bfi_{16} P$. In the same way one gets that
$p \bfi_{16} P$ implies $p \bfi_V P$. This proves the claim.
\end{proof}
\end{Num}

\begin{Thm} \label{la.TheFinalTheorem}
Let $\calh$ be the set of strongly isotropic points. Then the relation
$\bot \subseteq \calh \times \calh$ defined by
$$\cc X \bot \cc Y \Longleftrightarrow X \times Y = 0$$
is a collinearity relation, and the incidence geometry $\Gamma_\bot$ associated
to $\bot$ is isomorphic to the real $E_6$ quadrangle $Q(E_6,\rr)$.
\proof
The point set $\calp$ of $Q(E_6,\rr)$ consists of all points
$\cc X$ ($X$ is a Veronese vector) which are incident with $P(H(X) \times
\herm_3^\cc)$. Hence the set $\calp$ equals the set $\calh$ of strongly
isotropic points (compare \ref{la.Veldkamp}). By
\ref{la.CollinearityE6Quadrangle} the collinearity relation $\bot$ in
$\Gamma$ describes the collinearity in $Q(E_6,\rr)$. Using
\ref{la.CollinearityRelation} we derive the claim.
\qed
\end{Thm}

\appendix

\section{Equations for strongly isotropic points} \label{sec5}

In this appendix we list concrete equations that describe the Veronese embedding of the $E_6$ quadrangle given in Theorem \ref{la.TheFinalTheorem}.

\begin{Prop}
Suppose $p = \cc X \in \calv$ is a point where $X=\hv{\xi}{x}$ is a Veronese
vector. Then $p$ is strongly isotropic if and only if the equations
\begin{eqnarray*}
  |\xi_1|^2 + |\xi_2|^2 + |\xi_3|^2 & = & \langle x_1 | \ovl{x}_1 \rangle
    + \langle x_2 | \ovl{x}_2 \rangle - \langle x_3 | \ovl{x}_3 \rangle \\
  |\xi_j|^2 + |\xi_k|^2 + \sigma_i \langle x_i | \ovl{x}_i \rangle & = & |\xi_i|^2 \\
  \ovl{\xi}_kx_j + \sigma_i\conj{{(x_k\ovl{x}_i)}} & = & - \sigma_j\xi_i\ovl{x}_j \\
  \ovl{\xi}_jx_k + \sigma_i\conj{{(\ovl{x}_ix_j)}} & = & - \sigma_k\xi_i\ovl{x}_k \\
  \sigma_i\xi_k\ovl{x}_i + \sigma_k\conj{{(x_j\ovl{x}_k)}} & = & - \ovl{\xi}_jx_i \\
  \sigma_i\xi_j\ovl{x}_i + \sigma_j\conj{{(\ovl{x}_jx_k)}} & = &- \ovl{\xi}_kx_i \\
  \sigma_j(tx_j)\conj{\ovl{x}_j} + \sigma_k\conj{\ovl{x}_k}(x_kt) + |\xi_i|^2t
    + \sigma_i \langle x_i | t \rangle \ovl{x}_i & = & \sigma_i \langle \ovl{x}_i |
    t \rangle x_i \\
  \sigma_j(x_it)\conj{\ovl{x}_j} - \sigma_k\xi_j\conj{{(t\ovl{x}_k)}}
    - \ovl{\xi}_i\conj{{(tx_k)}} & = & \sigma_j \langle \ovl{x}_j | t \rangle x_i \\
  - \sigma_j\xi_k\conj{{(\ovl{x}_jt)}} + \sigma_k\conj{\ovl{x}_k}(tx_i)
    - \ovl{\xi}_i\conj{{(x_jt)}} & = & \sigma_k \langle \ovl{x}_k | t \rangle x_i
\end{eqnarray*}
hold for all $t \in \oo^\cc$ and $(i,j,k) \in \{ (1,2,3), (2,3,1), (3,1,2) \}$.
\proof
Let $T=\hv{\tau}{t} \in
\herm_3^\cc$ and put $Y = \hv{\eta}{y} = 2X \times T$. Then we
obtain from \ref{la.CrossProduct}
\begin{eqnarray*}
  \eta_i & = & \xi_j\tau_k + \xi_k\tau_j - \langle x_i|t_i \rangle \\
  y_i & = & \conj{{(x_jt_k)}} + \conj{{(t_jx_k)}} - \xi_it_i - \tau_ix_i
\end{eqnarray*}
for all $(i,j,k) \in \{ (1,2,3), (2,3,1), (3,1,2) \}$. Define
$(\sigma_1,\sigma_2,\sigma_3) = (-1,-1,1)$. Then we have
$$H(X) = (\ovl{\xi}_1, \ovl{\xi}_2, \ovl{\xi_3}; \sigma_1\ovl{x}_1,
\sigma_2\ovl{x}_2, \sigma_3\ovl{x}_3).$$
If we put
$$Z = 2H(X) \times Y = 4H(X) \times (X \times T)$$
where $Z=\hv{\zeta}{z}$, then we compute
\begin{eqnarray*}
  \zeta_i & = & \ovl{\xi}_j\eta_k + \ovl{\xi}_k\eta_j
    - \langle \sigma_i\ovl{x}_i|y_i \rangle \\
  & = & \ovl{\xi}_j  ( \xi_i\tau_j + \xi_j\tau_i - \langle x_k|t_k \rangle )
    + \ovl{\xi}_k ( \xi_k\tau_i + \xi_i\tau_k - \langle x_j|t_j \rangle ) \\
  & & - \langle \sigma_i\ovl{x}_i | \conj{{(x_jt_k)}} + \conj{{(t_jx_k)}} - \xi_it_i
    - \tau_ix_i \rangle \\
  & = & (\ovl{\xi}_j\xi_j + \ovl{\xi}_k\xi_k + \langle \sigma_i\ovl{x}_i |
    x_i \rangle ) \tau_i + \ovl{\xi}_j\xi_i\tau_j + \ovl{\xi}_k\xi_i\tau_k \\
  & & + \langle \xi_i\sigma_i\ovl{x}_i | t_i \rangle - \langle \ovl{\xi}_kx_j
    + \conj{{(x_k\sigma_i\ovl{x}_i)}} | t_j \rangle - \langle \ovl{\xi}_jx_k
    + \conj{{(\sigma_i\ovl{x}_ix_j)}} | t_k \rangle
\end{eqnarray*}
and
\begin{eqnarray*}
  z_i & = & \conj{{(\sigma_j\ovl{x}_jy_k)}} + \conj{{(y_j\sigma_k\ovl{x}_k)}}
    - \ovl{\xi}_iy_i - \eta_i\sigma_i\ovl{x}_i \\
  & = & \conj{{ \Big( \sigma_j\ovl{x}_j \big( \conj{{(x_it_j)}} + \conj{{(t_ix_j)}}
    - \xi_kt_k - \tau_kx_k \big) \Big) }} \\
  & & + \conj{{ \Big( \big( \conj{{(x_kt_i)}} + \conj{{(t_kx_i)}} - \xi_jt_j - \tau_jx_j \big)
    \sigma_k\ovl{x}_k \Big) }} \\
  & & - \ovl{\xi}_i \big( \conj{{(x_jt_k)}} + \conj{{(t_jx_k)}} - \xi_it_i
    - \tau_ix_i \big) - \big( \xi_j\tau_k + \xi_k\tau_j - \langle x_i | t_i \rangle \big)
    \sigma_i\ovl{x}_i \\
  & = & \ovl{\xi}_i\tau_ix_i - \tau_j \big( \conj{{(x_j\sigma_k\ovl{x}_k)}}
    + \xi_k\sigma_i\ovl{x}_i \big) - \tau_k \big( \conj{{(\sigma_j\ovl{x}_jx_k)}}
    + \xi_j\sigma_i\ovl{x}_i \big) \\
  & & + (t_ix_j)\sigma_j\conj{\ovl{x}_j} + \sigma_k\conj{\ovl{x}_k}(x_kt_i)
    + \ovl{\xi}_i\xi_it_i + \langle x_i | t_i \rangle \sigma_i\ovl{x}_i \\
  & & + (x_it_j)\sigma_j\conj{\ovl{x}_j} - \xi_j\conj{{(t_j\sigma_k\ovl{x}_k)}}
    - \ovl{\xi}_i\conj{{(t_jx_k)}} \\
  & & - \xi_k\conj{{(\sigma_j\ovl{x}_jt_k)}} + \sigma_k\conj{\ovl{x}_k} (t_kx_i)
    - \ovl{\xi}_i\conj{{(x_jt_k)}}
\end{eqnarray*}
for all $(i,j,k) \in \{ (1,2,3), (2,3,1), (3,1,2) \}$. On the other hand, if
we put $Z = h(T,X)X$, then we have
$$\zeta_i = \tau_i\ovl{\xi}_i\xi_i + \tau_j\ovl{\xi}_j\xi_i +
\tau_k\ovl{\xi}_k\xi_i + \langle t_i | \sigma_i\ovl{x}_i \rangle \xi_i +
\langle t_j | \sigma_j\ovl{x}_j \rangle \xi_i + \langle t_k |
\sigma_k\ovl{x}_k \rangle \xi_i$$
and
$$z_i = \tau_i\ovl{\xi}_ix_i + \tau_j\ovl{\xi}_jx_i + \tau_k\ovl{\xi}_kx_i +
\langle t_i | \sigma_i\ovl{x}_i \rangle x_i + \langle t_j | \sigma_j\ovl{x}_j
\rangle x_i + \langle t_k | \sigma_k\ovl{x}_k \rangle x_i$$
for each $(i,j,k) \in \{ (1,2,3), (2,3,1), (3,1,2) \}$. Now a comparison of
both sides of the equation $4H(X) \times (X \times T) = h(T,X)X$ yields the
claim.
\qed
\end{Prop}

\bibliographystyle{num}

\bigskip\makeatletter
\raggedright
Ralf Gramlich\\
TU Darmstadt, FB Mathematik AG 5, Schlossgartenstr. 7, 64289 Darmstadt, Germany.\\
{\tt gramlich@mathematik.tu-darmstadt.de}\\
The University of Birmingham, School of Mathematics, Edgbaston,
Birmingham B15 2TT, United Kingdom.\\
{\tt ralfg@maths.bham.ac.uk}

\medskip
Linus Kramer\\
Mathematisches Institut, 
Universit\"at M\"unster,
Einsteinstr. 62,
48149 M\"unster,
Germany.\\
{\tt linus.kramer{@}math.uni-muenster.de}\\\smallskip

\end{document}

%% file: e6i.tex
\setlength{\unitlength}{3947sp}%
\begingroup\makeatletter\ifx\SetFigFont\undefined%
\gdef\SetFigFont#1#2#3#4#5{%
  \reset@font\fontsize{#1}{#2pt}%
  \fontfamily{#3}\fontseries{#4}\fontshape{#5}%
  \selectfont}%
\fi\endgroup%
\begin{picture}(2266,1080)(983,-538)
\thicklines
\put(1036,-241){\circle{90}}
\put(1576,-241){\circle{90}}
\put(2116,-241){\circle{90}}
\put(2656,-241){\circle{90}}
\put(2116,299){\circle{90}}
\put(3196,-241){\circle{90}}
\put(1081,-241){\line( 1, 0){450}}
\put(2071,-241){\line(-1, 0){450}}
\put(2161,-241){\line( 1, 0){450}}
\put(2116,-196){\line( 0, 1){450}}
\put(3151,-241){\line(-1, 0){450}}
\put(991,-511){\makebox(0,0)[lb]{\smash{\SetFigFont{12}{14.4}{\rmdefault}{\mddefault}{\updefault}1}}}
\put(2611,-511){\makebox(0,0)[lb]{\smash{\SetFigFont{12}{14.4}{\rmdefault}{\mddefault}{\updefault}5}}}
\put(1531,-511){\makebox(0,0)[lb]{\smash{\SetFigFont{12}{14.4}{\rmdefault}{\mddefault}{\updefault}3}}}
\put(2071,-511){\makebox(0,0)[lb]{\smash{\SetFigFont{12}{14.4}{\rmdefault}{\mddefault}{\updefault}4}}}
\put(3151,-511){\makebox(0,0)[lb]{\smash{\SetFigFont{12}{14.4}{\rmdefault}{\mddefault}{\updefault}6}}}
\put(2116,434){\makebox(0,0)[lb]{\smash{\SetFigFont{12}{14.4}{\rmdefault}{\mddefault}{\updefault}2}}}
\end{picture}

%% file: e6ii.tex
\setlength{\unitlength}{3947sp}%
\begingroup\makeatletter\ifx\SetFigFont\undefined%
\gdef\SetFigFont#1#2#3#4#5{%
  \reset@font\fontsize{#1}{#2pt}%
  \fontfamily{#3}\fontseries{#4}\fontshape{#5}%
  \selectfont}%
\fi\endgroup%
\begin{picture}(1936,1575)(2070,-943)
\thicklines
\put(2251,389){\circle{90}}
\put(3331,-106){\circle{90}}
\put(3871,-106){\circle{90}}
\put(2791,-646){\circle{90}}
\put(2251,-646){\circle{90}}
\put(2789,386){\circle{94}}
\put(3871,-106){\oval(180,630)}
\put(2258,-128){\oval(360,1260)}
\put(2296,389){\line( 1, 0){450}}
\put(2836,389){\line( 1,-1){450}}
\put(3286,-151){\line(-1,-1){450}}
\put(3376,-106){\line( 1, 0){450}}
\put(2746,-646){\line(-1, 0){450}}
\put(2746,-916){\makebox(0,0)[lb]{\smash{\SetFigFont{12}{14.4}{\rmdefault}{\mddefault}{\updefault}5}}}
\put(3376,-16){\makebox(0,0)[lb]{\smash{\SetFigFont{12}{14.4}{\rmdefault}{\mddefault}{\updefault}4}}}
\put(2746,524){\makebox(0,0)[lb]{\smash{\SetFigFont{12}{14.4}{\rmdefault}{\mddefault}{\updefault}3}}}
\put(2206,-916){\makebox(0,0)[lb]{\smash{\SetFigFont{12}{14.4}{\rmdefault}{\mddefault}{\updefault}6}}}
\put(2251,524){\makebox(0,0)[lb]{\smash{\SetFigFont{12}{14.4}{\rmdefault}{\mddefault}{\updefault}1}}}
\put(4006,-196){\makebox(0,0)[lb]{\smash{\SetFigFont{12}{14.4}{\rmdefault}{\mddefault}{\updefault}2}}}
\end{picture}

%% file: e6.bbl
\begin{thebibliography}{99}%{SBG{\etalchar{+}}95}

\bibitem{Adams}
J. Frank Adams.
\newblock {\em Lectures on exceptional Lie groups.}
\newblock Chicago Lecture Notes in Mathematics,
University of Chicago Press, 1996.

\bibitem{DiagramGeometries}
Francis Buekenhout and Arjeh~M. Cohen.
\newblock {\em Diagram geometry.}
\newblock Book in preparation, a perliminary version is available at
A. Cohen's web site,
{\tt http://www.win.tue.nl/$\sim$amc/buek/}

\bibitem{JordanAlgebren}
Hel Braun and Max Koecher.
\newblock {\em Jordan-Algebren}.
\newblock Springer-Verlag, Berlin - Heidelberg - New York, 1966.

\bibitem{PointLineSpaces}
Arjeh~M. Cohen.
\newblock Point-line spaces related to buildings.
\newblock In Francis Buekenhout, editor, {\em Handbook of Incidence Geometry}.
  Elsevier, Amsterdam - Lausanne - New York, 1995.

\bibitem{ProjektiveOktavenebene}
Hans Freudenthal.
\newblock Oktaven, {A}usnahmegruppen und {O}ktavengeometrie.
\newblock {\em Geom. Dedicata}, 19:7--63, 1985.

\bibitem{Greub}
Werner H. Greub.
\newblock {\em Multilinear Algebra}.
\newblock Springer-Verlag, Berlin - Heidelberg - New York, 2nd edition, 1978.

\bibitem{PrinciplesOfAlgebraicGeometry}
Phillip Griffiths and Joseph Harris.
\newblock {\em Principles of Algebraic Geometry}.
\newblock John Wiley \& Sons, New York - Chichester - Brishane, 1978.

\bibitem{Helgason}
Sigurdur Helgason.
\newblock {\em Differential Geometry, Lie Groups, and Symmetric Spaces}.
\newblock Academic Press, New York - London, 1978.

\bibitem{ClassicalGroups}
Alexander~J. Hahn and O.~Timothy O'Meara.
\newblock {\em The Classical Groups and K-Theory}.
\newblock Springer-Verlag, Berlin - Heidelberg - New York, 1989.

%\bibitem{CompositionAlgebrasAutomorphisms}
%Nathan Jacobson.
%\newblock Composition algebras and their automorphisms.
%\newblock {\em Rend. Circ. Mat. Palermo (2)}, 7:55--80, 1958.

\bibitem{BasicAlgebraI}
Nathan Jacobson.
\newblock {\em Basic Algebra I}.
\newblock W. H. Freeman and Company, New York, 2nd edition, 1985.

%\bibitem{BasicAlgebraII}
%Nathan Jacobson.
%\newblock {\em Basic Algebra II}.
%\newblock W. H. Freeman and Company, New York, 2nd edition, 1989.

\bibitem{Kurth:2000}
Torsten Kurth.
\newblock {\em On a real form of $E_6$ and its related generalized quadrangle}.
\newblock Diplomarbeit, W\"urzburg 2000.

\bibitem{LangAlgebra}
Serge Lang.
\newblock {\em Algebra}.
\newblock Addison-Wesley, Reading, Massachusetts, 3rd edition, 1993.

\bibitem{PolarSpacesOfRank3}
Bernhard M{\"u}hlherr.
\newblock A geometric approach to non-embeddable polar spaces of rank $3$.
\newblock {\em Bull. Soc. Math. Belg.}, XLII:577--594, 1990.

\bibitem{MHerr}
Bernhard M{\"u}hlherr.
\newblock {\em Some Contributions to the Theory of Buildings Based on the Gate
  Property}.
\newblock PhD thesis, Universit\"at T\"ubingen, 1994.

\bibitem{Payne/Thas:1984}
Stanley Payne, Joseph A.\ Thas.
\newblock {\em Finite generalized quadrangles}.
\newblock  Pitman, Boston 1984.

\bibitem{CliffordAlgebras}
Ian~R. Porteous.
\newblock {\em Clifford Algebras and the Classical Groups}.
\newblock Cambridge University Press, Cambridge, 1995.

\bibitem{CompactProjectivePlanes}
Helmut Salzmann, Dieter Betten, Theo Grundh\"ofer, Hermann H\"ahl, Rainer
  L\"owen, and Markus Stroppel.
\newblock {\em Compact Projective Planes. With an Introduction to Octonion
  Geometry.}
\newblock Walter de Gruyter, Berlin, 1995.

\bibitem{BasicAlgebraicGeometryI}
Igor~R. Shafarevich.
\newblock {\em Basic Algebraic Geometry 1}.
\newblock Springer-Verlag, Berlin - Heidelberg - New York, 2nd edition, 1994.

\bibitem{ClassOfJordanAlgebras}
Tonny~A. Springer.
\newblock On a class of {J}ordan algebras.
\newblock {\em Indag. Math.}, 21:254--264, 1959.

\bibitem{ClassOfCubicForms}
Tonny~A. Springer.
\newblock Characterization of a class of cubic forms.
\newblock {\em Indag. Math.}, 24:259--265, 1962.

\bibitem{GeometricAlgebraProjectiveOctavePlanes}
Tonny~A. Springer.
\newblock On the geometric algebra of the octave planes.
\newblock {\em Indag. Math.}, 24:451--468, 1962.

\bibitem{HjelmslevI}
Tonny~A. Springer and Ferdinand~D. Veldkamp.
\newblock On {H}jelmslev-{M}oufang planes.
\newblock {\em Math. Z.}, 107:249--263, 1968.

\bibitem{GeometryOfClassicalGroups}
Donald~E. Taylor.
\newblock {\em The Geometry of the Classical Groups}.
\newblock Heldermann Verlag, Berlin, 1992.

%\bibitem{Tabellen}
%Jacques Tits.
%\newblock {\em Tabellen zu den einfachen Lie Gruppen und ihren Darstellungen},
%  volume~40 of {\em Lecture Notes in Mathematics}.
%\newblock Springer-Verlag, Berlin - Heidelberg - New York, 1967.

\bibitem{BuildingsOfSphericalType}
Jacques Tits.
\newblock {\em Buildings of Spherical Type and Finite BN-pairs}, volume 386 of
  {\em Lecture Notes in Mathematics}.
\newblock Springer-Verlag, Berlin - Heidelberg - New York, 1974.

\bibitem{TitsWeiss}
Jacques Tits\ and\ Richard M.\ Weiss, {\it Moufang polygons}, Springer, Berlin, 2002.% MR1938841 (2003m:51008)

\bibitem{HjelmslevIII}
Ferdinand~D. Veldkamp.
\newblock Unitary groups in {H}jelmslev-{M}oufang planes.
\newblock {\em Math. Z.}, 108:288--312, 1969.

\bibitem{GeneralizedPolygons}
Hendrik Van~Maldeghem.
\newblock {\em Generalized Polygons}.
\newblock Birkh\"auser, Basel - Boston - Berlin, 1998.

\bibitem{WeissE6}
Richard M.\ Weiss, Moufang quadrangles of type $E\sb 6$ and $E\sb 7$, J. Reine Angew. Math. {\bf 590} (2006), 189--226.% MR2208133 (2006j:51005)

\bibitem{Weiss:2006}
Richard M.\ Weiss, {\em Quadrangular algebras}, Princeton University Press, Princeton, 2006.

\end{thebibliography}
